%% file: abkm2020.tex
\documentclass[11pt]{article}

\usepackage{amsmath}
\usepackage{amsfonts}
\usepackage{hyperref}
\usepackage{enumitem}
\usepackage{natbib}
\usepackage{graphicx}
\usepackage{multirow}

\DeclareMathOperator*{\argmin}{\mathrm{argmin}}

\begin{document}

\title{Relax-and-fix heuristics applied to a real-world lot-sizing and scheduling problem in the personal care consumer goods industry\thanks{This work has been partially supported by the Brazilian agencies FAPESP (grants 2013/07375-0, 2016/01860-1, and 2018/24293-0) and CNPq (grants 306083/2016-7 and 302682/2019-8).}}

\author{
    K. A. G. Araujo\thanks{Department of Applied Mathematics, Institute of
    Mathematics and Statistics, University of S\~ao Paulo, Rua do
    Mat\~ao, 1010, Cidade Universit\'aria, 05508-090, S\~ao Paulo, SP,
    Brazil. e-mail: kennedy94@ime.usp.br}
 \and 
    E. G. Birgin\thanks{Department of Computer Science, Institute of
    Mathematics and Statistics, University of S\~ao Paulo, Rua do
    Mat\~ao, 1010, Cidade Universit\'aria, 05508-090, S\~ao Paulo, SP,
    Brazil. e-mail: egbirgin@ime.usp.br}
 \and
    M. S. Kawamura\thanks{Department of Production Engineering, Polytechnic 
    School, University of S\~ao Paulo, Av. Prof. Luciano Gualberto, 1380, 
    Cidade Universit\'aria, 05508-010 São Paulo, SP, Brazil. e-mails: 
    marcio.kawamura@poli.usp.br, dronconi@usp.br}
 \and
    D. P. Ronconi\footnotemark[4]
}

\date{July 22, 2021}

\maketitle

\begin{abstract}
This paper addresses an integrated lot-sizing and scheduling problem in the industry of consumer goods for personal care, a very competitive market in which the good customer service level and the cost management show up in the competition for the clients. In this research, a complex operational environment composed of unrelated parallel machines with limited production capacity and sequence-dependent setup times and costs is studied. There is also a limited finished-goods storage capacity, a characteristic not found in the literature. Backordering is allowed but it is extremely undesirable. The problem is described through a mixed integer linear programming formulation. Since the problem is NP-hard, relax-and-fix heuristics with hybrid partitioning strategies are investigated. Computational experiments with randomly generated and also with real-world instances are presented. The results show the efficacy and efficiency of the proposed approaches. Compared to current solutions used by the company, the best proposed strategies yield results with substantially lower costs, primarily from the reduction in inventory levels and better allocation of production batches on the machines.\\

\noindent
\textbf{Keywords:} lot sizing and scheduling, mixed integer linear programming models, relax-and-fix, real-world instances.\\
\end{abstract}

\section{Introduction}

According to the Brazilian Toiletry, Perfumery and Cosmetic Association (ABIHPEC), the Brazilian market for personal care, perfumery, and cosmetics showed a growth of $2.2\%$ in 2020 above the $-4.5\%$ of the overall industry and $-4.1\%$ of the GDP (Gross Domestic Product) (both heavily affected by the COVID-19 pandemic); and a CAGR (Compound Annual Growth Rate) of $1.7\%$ over the last 10 years compared to a CAGR of $-2.1\%$ for industry overall and $-0.3\%$ for GDP over the same period~\citep{abihpec}. Among the factors that have contributed to this accelerated growth are the growing participation of women in the labor market, the frequent releases of new products, the higher productivity by using cutting-edge technology, and, more recently, the essentiality of the segment's products in the combat of COVID-19 pandemic. This segment is comprised of a large number of small manufacturers and some large companies. In 2020, Brazil held the fourth biggest consumer market of this sector in the world with US\$~$23.7$ billions. Inside this market, the segment studied in this work is the disposable personal care one, composed of diapers, sanitary pads, toilet paper, towels, and tissues. These products are characterized by (i) frequent consumption by the population, (ii) relatively low cost, and (iii) very similar quality among different brands. Therefore the process of production planning, scheduling and control performs an important role in the companies to guarantee productivity, cost control, and adequate customer service level. In most companies, the main responsibility of this activity is to simultaneously analyze several relevant information such as the seasonality of the demand, the seasonality of the supply of raw materials, perishability of the finished products and the raw material, and the manufacturing capacity to develop a production plan that optimizes the use of productive resources while meets the demand for manufactured products.

Typical production planning decisions as lot sizing and scheduling significantly influence the results of a company by maximizing the fulfillment of sales orders, adjusting correctly the inventory levels, and incurring in proper operating costs. The process of lot sizing consists of determining how much to produce of each product in each period in order to meet a projected demand under the existing conditions and operational capabilities. On the other hand, scheduling means determining when and in which sequence these lots have to be produced in order to maximize the productive resources efficiency and to meet the deadlines. Inefficiencies in these processes can cause overstocking of finished products, unfulfilled sales orders, loss of perishable material, not accomplished due dates, significant reduction in the productive capacity of the production line, stocks of finished products accumulated in advance, higher preparation machine costs, among others.

\begin{figure}[!ht]
\centering
\includegraphics[scale=1.0]{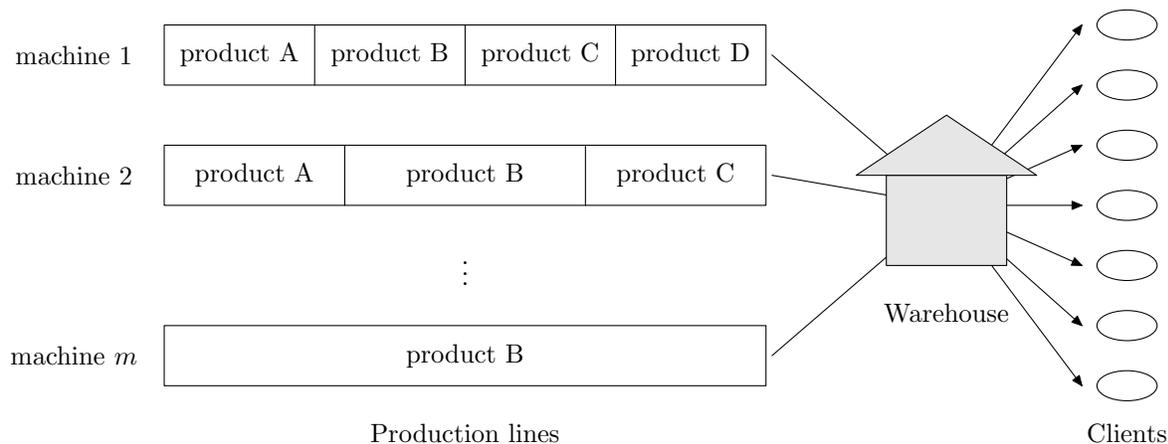}
\caption{Schematically representation of the production process.}
\label{fig1}
\end{figure}

In this work, a case study of a large company in the market of disposable personal care products is performed. The company has factories in the south and southeast of Brazil and its market is also geographically concentrated in these regions. Its main product categories are diapers, feminine sanitary pads, and toilet paper. These products are manufactured by the company itself; and their demands are predominantly affected by the own company marketing activities and by the competition. Some other manufactured products have a strongly seasonal demand as tissue paper, which are widely used in the winter due to higher incidence of respiratory diseases, and disposable napkins and paper towels, sold mostly on holiday periods like Christmas, Easter, and New Year, among others. The production configuration of the company consists of machines of different models that have been acquired over time, as demand was increasing. These machines have performance profiles distinct from each other with regard to efficiency, speed, and cost. Moreover,  by technical constraints,  not all products can be manufactured in all machines. Between the production of batches of different products on the same line, machine preparation (setup) is required, which may be simple and fast as changing the consuming packing or complex and slow if it involves change of raw material or reconfiguration of the parameters of the product, among others. Thus, these setups depend on the product to be produced and also on the product which was being produced in line in the previous batch. The products manufactured are shipped to a limited capacity distribution center located close to the factories. Figure~\ref{fig1} schematically represents the process.

In the literature, the integrated lot sizing and scheduling problem with dependent setup times and costs on a single machine (or a singe line) is called GLSPST (General Lot Sizing and Scheduling Problem with sequence-dependent Setup Time); while the same problem involving parallel machines is called GLSPPL (General Lot Sizing and Scheduling Problem with Parallel Machines). The GLSP with non-zero minimum lotsizes and without setup times has been proved to be NP-complete by \cite{fleischmann1997general}. Consequently, the GLSP and the GLSPPL with sequence setup times are also NP-hard~\citep{meyr2000simultaneous,meyr2002simultaneous}. 
The focused problem  can be seen as the GLSPPL with sequence-dependent setup times and costs, non-identical parallel machines, in addition to specific characteristics of the environment addressed, such as limited warehousing capacity for finished products, machine eligibility constraints,  backorder costs, and setup state conservation among adjacent periods. Up to the authors acknowledge, few studies addressed the GLSPPL with all the characteristics of the focused real-word problem.

In this research, a heuristic procedure able to tackle real-world instances of a large company in the market of hygiene products  is  presented. The proposed solution method is a customized relax-and-fix heuristic that decompose the original mixed integer linear programming (MILP) model of the problem into less complex subproblems that are successively solved. This strategy was selected due to the fact that this heuristic method had been successfully applied to similar problems (see for example \citet{clark2000rolling, beraldi2008rolling, lang2011fix, soler2020}). Furthermore, this MILP-based approach is mentioned by \cite{copil} as the next generation of solution methods for the MILP models of the lot sizing and scheduling problems. Eleven different solving approaches are presented. Nine of them are problem-dependent strategies that divide the main problem considering several metrics associated with machines, products, and periods; and two of them are problem-independent strategies. The aim of the study is not only to solve a complex real-world problem but also to assess in which way the form of partitioning the problem into simpler subproblems and the sequence in which subproblems are solved affect the performance of the relax-and-fix heuristic.

The rest of this paper is organized as follows. A literature review is presented in Section~\ref{seclierature}. In Section~\ref{secmodel}, a  MILP formulation of the problem is given. In Section~\ref{secrelax} the considered relax-and-fix heuristics are described. In Section~\ref{secnum} numerical experiments with randomly generated instances and real-world instances of the company are conducted.
The proposed heuristic methods are compared against a commercial solver applied to the MILP model presented in Section~\ref{secmodel}; while, for the real-world instances, a comparison with the solutions adopted by the company is also presented. Conclusions are given in the last section.

\section{Literature review}
\label{seclierature}

Due to the relevance of the production planning for the industry and its complexity, many authors have addressed the lot sizing and scheduling problem with different characteristics. See classical reviews about the theme in \cite{drexl1997lot}, \cite{karimi2003capacitated}, \cite{zhu2006scheduling}, \cite{quadt2008capacitated}, \cite{robinson2009coordinated}, and for a recent review see \cite{copil}. A brief description of the works that address that GLSPPL with specific characteristics is presented below. Differences with the problem considered in this work are highlighted.

\cite{kang1999lotsizing} treat the GLSPPL considering the minimization of setup and inventory costs using column generation techniques. In this work the cost of sequence-dependent setups is considered, but setup times are ignored. \cite{clark2000rolling} studies the GLSPPL whose objective is the minimization of storage and backordering costs, but without considering setup costs. The proposed mathematical model is based on the premise that the maximum number of setups per period is predetermined. A rolling-horizon method and relax-and-fix heuristics are applied. However, the presented computational results shows that only small problems could be solved in reasonable time. Aiming to minimize production, inventory and setup costs without allowing delivery delays, \cite{meyr2002simultaneous} considers the metaheuristics Threshold Accepting and Simulated Annealing in small real instances with identical machines. \cite{beraldi2008rolling} develop rolling-horizon and relax-and-fix heuristics to solve the GLSPPL in the textile and fiberglass industry environment. Unlike the mathematical model developed by \cite{meyr2002simultaneous}, the authors introduce a compact formulation for the case of identical machines considering the  setup costs, but neglecting the setup times. 

\cite{jozefowska2008optimization} develop a decision support system applied to a Polish company that manufacture plastic pipes, whose production environment is composed of unrelated parallel machines. This system is based on a multi-objective model which includes, among its criteria, the maximization of the machine utilization and the minimization of the deviation between the production schedule and the S\&OP (Sales and Operations Plan) and the amount of products below the required level of safety stock. The model is solved using a genetic algorithm after having its solution space reduced by adding restrictions suggested by experienced planners. \cite{mateus2010capacitated} approach the lot sizing and scheduling problem considering a factory of refractory bricks with different machines in parallel.  Unlike the problem addressed in this work, the authors do not consider setup carry-over, i.e.\ the preparation of the machines is not maintained from one period to the next. The proposed iterative solution method is composed by two modules: the first solves exactly the problem of lot sizing considering aggregate capacity and estimated setup times; while the second searches a feasible sequencing for pre-sized lots through a GRASP meta-heuristic. \cite{meyr2013decomposition} deal with the GLSPPL composed by heterogeneous parallel machines with the objective of minimizing inventory and sequence-dependent setup and production costs, without backlogging. The authors propose an heuristic which iteratively decomposes the multi-line problem into a series of single-line subproblems, which can be easily solved  by the heuristic TADR proposed by \cite{meyr2000simultaneous}. Two strategies of decomposing the problem into subproblems are proposed: (i) a strategy based on priority rules and (ii) another based on aggregating the original problem, solving the aggregate problem, and disaggregating the results to define the demand and initial inventory for each line. \cite{xiao2015hybrid} also examine the parallel-machine lot-sizing and scheduling problem with sequence-dependent setup times.  In addition, release times of the items and  machine eligibility and preference constraints are considered. However, the setup cost is sequence-independent and the production costs are not considered. The authors initially propose a MILP model; next  the original problem is decomposed into a lot-sizing subproblem and a set of single-machine scheduling subproblems by Lagrangian decomposition. A Lagrangian-based heuristic algorithm, which incorporates the simulated annealing algorithm to improve the solution of the scheduling problem, is proposed. 

Considering the GLSPPL with secondary resources, \cite{dastidar2005scheduling} approach a problem on the health-related injection molding industry where a number of resources, available in limited quantities (e.g. grinders and driers), are installed on machines to perform certain tasks. The authors tackle the problem through a two-stage decomposition. In the first phase, the machines are grouped according to their criticality (machines that produce products with fewer machine options to be processed are more critical). In the second phase, the sub-problems represented by a MILP model, are solved using an open-source solver. \cite{alme2011} present MILP models to tackle the synchronization of a secondary resource in lot-sizing and scheduling problems with parallel unrelated machines. The machines have to be equipped with a special kind of resource (e.g. a tool) with limited capacity and multiple products can be produced with the same tool, therefore the tool changeovers costs and times are considered instead of product changeovers. A GLSPPL inspired by a real-world production environment in the food industry is addressed by \cite{soler2020}. In the problem at hand, due to the scarcity of resources, only a subset of production lines can operate simultaneously and these lines need to be assembled at each production period. Due to this feature, no setup carryover between adjacent periods exist, only a subset of products can be produced in a given period, and for each product, only one production line is capable of producing it. The authors define a branching rule to accelerate the performance of the branch-and-bound algorithm of the CPLEX solver and a relax-and-fix procedure is also implemented. More recently, \cite{yasmin2020} tackle a GLSPPL to maximize the profit of assembled products  where pieces are produced using auxiliary equipment (molds) to form finished products. Each piece may be processed in a set of molds with different production rates on various machines. An iterative heuristic which decouples the lot-sizing and scheduling decisions is proposed.  First, the lot-size of the products is determined through the solution of a MILP model where a mold can be used in more than one machine at a time. Next, another MILP model is presented that, based on the lot-size and mold-machine assignment previously obtained, allows to determine the molds schedule in the machines in each planning period. \cite{dearmas2020} address the GLSPPL in a context of pipe insulation manufacturing in which a limited number of resources must be shared by parallel machines and sets of make-to-stock (MTS) and make-to-order (MTO) units must be considered simultaneously. The stock level of each MTS stock keeping unit (SKU) must be within specified limits at the end of each period, while each MTO SKU must be produced within its specified time frame. The objective is to maximize the total amount of production in a given planning horizon considering sequence-dependent setup times; the setup costs are not considered. The authors propose a two-phase procedure. First, a MIP model in which the setup times are considered independent of the sequence is formulated and solved. Next, a one-pass heuristic is applied to search for setup-time savings by reordering the products assigned by the MIP model to each machine in each period. An overview on simultaneous lot sizing and scheduling involving secondary resources can be found in \cite{worbe2019}.

\section{Mathematical Model} \label{secmodel}

In the GLSPPL, $n$ types of products are manufactured in a shop floor composed of~$m$ different machines over a horizon divided into~$T$ time periods. In a period $t$ (the time interval that goes from instant $t-1$ to instant~$t$), a machine can produce more than one type of product, provided its time availability is not exceed. Machines have different production rates and efficiency levels. Some machines, due to their manufacturer, model, or preservation status can achieve a higher production speed with lower scrap rates than others. Consequently, the production time~$p_{i \ell}$ and the production cost $c^P_{i \ell}$ to produce product~$i$ on machine~$\ell$ depend on the product and the machine. Furthermore, by technical constraints, not all products can be manufactured in all machines. Between production batches of different types of products, it is necessary a time to prepare the machine and to set the correct product parameters, which generally generates loss of material. These setup times $e_{i j \ell}$ and, consequently, their respective costs $c^S_{ij\ell}$, are sequence dependent. Preemption is not permitted. The demand~$d_{it}$ of each product~$i$ at the end of period~$t$ is dynamic and deterministic, i.e.\ it is known and it varies along the time. Backordering is allowed but each backordered unit of product~$i$ is penalized by~$g_i$ per period of delay. Finished goods are transferred from the factory to a centralized distribution center with limited warehousing capacity~$C^W$. Each stored unit of product~$i$ costs an inventory cost~$h_i$ per period. In the considered problem, besides defining the quantities of each product to be produced and the production sequence, the solution determines in which machine each batch is manufactured in order to minimize the sum of the costs of inventory, setup, production, and backordering.

The MILP formulation presented in the current section is based on the MILP formulation introduced in~\citet{meyr2002simultaneous}, that does not consider the allowance of backorders, machine eligibility constraints, and the limited capacity of warehousing. In the formulation, many products can be produced in a machine in a certain period of time. The machine availability is consumed by the time necessary to setup the machine and to produce the lot. The setup time and cost depend on the sequence of products that are produced in the machine. The planning horizon is divided in $T$ periods. Within each period $t$,  each machine $\ell$ has $w_{\ell t}$ variable-length subperiods. A machine can produce a single type of product within each subperiod; and the duration of the subperiod is given by the duration of the production of the lot that is being produced. There may be subperiods with zero length, even though the machine is set up to produce some product. The division into subperiods determines the sequence of the jobs on each machine and also defines the associated setup times and costs.

In order to present the proposed model we introduce the following notation:

\noindent
\textbf{Main constants:}
\begin{description}[topsep=2pt,itemsep=1pt]
\item[]$m$: number of machines,
\item[]$n$: number of products,
\item[]$T$:	time horizon (assumed to start at~$0$).
\end{description}

\noindent
\textbf{Indexes:}
\begin{description}[topsep=2pt,itemsep=1pt]
\item[]$i$, $j$: products,
\item[]$\ell$: machines,
\item[]$t$:	periods (a period $t$ corresponds to the time interval between instants $t-1$ and $t$),
\item[]$s$:	subperiods.
\end{description}

\noindent
\textbf{Sets:}
\begin{description}[topsep=2pt,itemsep=1pt]
\item[]${\cal I} = \{ 1, \dots, n\}$: products,
\item[]${\cal L} = \{ 1, \dots, m\}$: machines,
\item[]${\cal I}_{\ell} \subseteq {\cal I}$: products that can be produced by machine $\ell$ ($\ell \in {\cal L}$),
\item[]${\cal L}_i \subseteq {\cal L}$: machines that can produce product $i$ ($i \in {\cal I}$),
\item[]${\cal T} = \{ 1, \dots, T\}$: periods,
\item[]${\cal S}_{\ell} = \{ 1, \dots, W_{\ell} \}$: subperiods of machine $\ell$ ($\ell \in {\cal L}$), where $w_{\ell t}$ is the number of subperiods of machine~$\ell$ within period $t$ ($\ell \in {\cal L}$, $t \in {\cal T}$) and $W_{\ell} = \sum_{t \in {\cal T}} w_{\ell t}$ ($\ell \in {\cal L}$),
\item[]${\cal S}_{\ell t} = \{ s \in S_{\ell} \; | \; \bar w_{\ell t} + 1 \leq s \leq \bar w_{\ell t} + w_{\ell t}\}$: subperiods of machine $\ell$ within period $t$ ($\ell \in {\cal L}$, $t \in {\cal T}$), where $\bar w_{\ell t} = 0$ for $t=1$ and $\bar w_{\ell t} = \bar w_{\ell,t-1} + w_{\ell,t-1}$ for $t = 2,\dots,T$.
\end{description}

\noindent
\textbf{Parameters:}
\begin{description}[topsep=2pt,itemsep=1pt]
\item[]$C^W$: capacity of warehousing,
\item[]$C^P_{\ell t}$: amount of time machine $\ell$ is available (for production and setup) within period $t$ ($\ell \in {\cal L}, t \in {\cal T}$),
\item[]$d_{it}$: demand of product $i$ at the end of period $t$, i.e.\ at instant $t$ ($i \in {\cal I}$, $t \in {\cal T}$),
\item[]$h_i$: inventory cost per period of a unit of product $i$ ($i \in {\cal I}$),
\item[]$g_i$: backordering cost per period of a unit of product $i$ ($i \in {\cal I}$),
\item[]$p_{i\ell }$: time required to produce a unit of product $i$ in machine $\ell$ ($\ell \in {\cal L}$, $i \in {\cal I_{\ell}}$),
\item[]$c^P_{i\ell }$: cost of producing a unit of product $i$ in machine $\ell$ ($\ell \in {\cal L}$, $i \in {\cal I_{\ell}}$),
\item[]$q^{\mathrm{lb}}_{i\ell }$: minimum lot of product $i$ that can be produced in machine $\ell$ ($\ell \in {\cal L}$, $i \in {\cal I_{\ell}}$),
\item[]$e_{ij\ell }$: setup time required to produce product $j$ immediately after $i$ in machine~$\ell$ ($\ell \in {\cal L}$, $i, j \in {\cal I_{\ell}}$),
\item[]$c^S_{ij\ell}$: cost of the setup required to produce product $j$ immediately after $i$ in machine~$\ell$ ($\ell \in {\cal L}$, $i, j \in {\cal I_{\ell}}$),
\item[]$I_{i0}^+$: inventory (quantity) of product $i$ at instant $t=0$ ($i \in {\cal I}$); it must satisfy $I_{i0}^+ \leq CW$,
\item[]$I_{i0}^-$: backordering (quantity) of product $i$ at instant $t=0$ ($i \in {\cal I}$),
\item[]$x_{i\ell 0}$: 1, if machine $\ell $ is prepared to produce product $i$ at instant $t=0$; 0, otherwise ($\ell \in {\cal L}$, $i \in {\cal I_{\ell}}$).
\end{description}

\noindent
\textbf{Variables:}
\begin{description}[topsep=2pt,itemsep=1pt]
\item[]$q_{i\ell s}$: quantity of product $i$ produced in machine $\ell$ within subperiod $s$ ($\ell \in {\cal L}$, $s \in S_{\ell}$, $i \in {\cal I_{\ell}}$), 
\item[]$x_{i\ell s}$: 1, if machine $\ell$ is prepared to produce product $i$ at the beginning of subperiod~$s$ (i.e.\ at instant $s-1$); 0, otherwise ($\ell \in {\cal L}$, $s \in S_{\ell}$, $i \in {\cal I_{\ell}}$),
\item[]$y_{ij\ell s}$:	1, if the setup required to produce product $j$ immediately after $i$ in machine $\ell$ occurs within subperiod $s$, i.e.\ between instants $s-1$ and $s$; 0, otherwise ($\ell \in {\cal L}$, $s \in S_{\ell}$, $i, j \in {\cal I_{\ell}}$),
\item[]$I_{it}^+$: inventory (quantity) of product $i$ at the end of period $t$, i.e.\ at instant $t$ ($i \in {\cal I}$, $t \in {\cal T}$),
\item[]$I_{it}^-$: backordering (quantity) of product $i$ at the end of period $t$, i.e.\ at instant $t$ ($i \in {\cal I}$, $t \in {\cal T}$).
\end{description}

The proposed mathematical formulation is presented below:
 
\begin{align}
\mbox{Minimize }\nonumber\quad & \\
\sum_{i \in {\cal I}} \sum_{t \in {\cal T}} h_i \, I_{it}^+
+ \sum_{i \in {\cal I}} \sum_{t \in {\cal T}} g_i \, I^-_{it} 
+ \sum_{\ell \in {\cal L}} \sum_{i \in {\cal I_{\ell}}} \sum_{j \in {\cal I_{\ell}}}  & \sum_{s \in S_{\ell}} c^S_{i j \ell}\, y_{i j \ell s}
+ \sum_{i \in {\cal I}} \sum_{\ell \in {\cal L}_i} \sum_{s \in S_{\ell}} c^P_{i \ell}\, q_{i\ell s} \label{fo}
\\
\mbox{subject to } \quad & \nonumber
\\
I^+_{i, t- 1} - I^-_{i, t-1} + \sum_{\ell  \in {\cal L}_i} \sum_{s \in S_{\ell t}} q_{i\ell s} - I^+_{it} + I^-_{it} = d_{it}, \quad & 
\quad i \in {\cal I}, \; t \in {\cal T}, \label{eq2}
\\ 
\sum_{i \in {\cal I}} I^+_{it} \leq C^W, \quad & 
\quad t \in {\cal T}, \label{eq3}
\\
\sum_{i \in {\cal I}_{\ell}}  \sum_{s \in S_{\ell t}} \left( p_{i\ell }\, q_{i\ell s} + \sum_{j \in {\cal I}_{\ell}} e_{ij\ell } \, y_{ij\ell s}\right) \leq C^P_{\ell t}, \quad & \quad \ell \in {\cal L}, \; t \in {\cal T}, \label{eq4}
\\
p_{i\ell } \, q_{i\ell s} \leq C^P_{\ell t} \, x_{i\ell s}, \quad & 
\quad \ell \in {\cal L}, \; t \in {\cal T}, \; s \in S_{\ell t}, \; i \in {\cal I}_{\ell} \label{eq5}
\\[2mm]
q_{i\ell s} \geq q^{\mathrm{lb}}_{i\ell }\, (x_{i\ell s} - x_{i,\ell ,s-1}), \quad & 
\quad \ell \in {\cal L}, \; t \in {\cal T}, \; s \in S_{\ell t}, \; i \in {\cal I}_{\ell} \label{eq6}
\\
\sum_{i \in {\cal I}_{\ell }} x_{i\ell s} = 1, \quad & 
\quad \ell \in {\cal L}, \; s \in S_{\ell}, \label{eq7}\\
y_{ij\ell s} \geq x_{i\ell ,s-1} + x_{j\ell s} - 1, \quad & \quad
\ell \in {\cal L}, \; s \in S_{\ell}, \; i,j \in {\cal I}_{\ell} \label{eq8}
\\
x_{i\ell s} \in \{0,1\}, \quad & \quad
\ell \in {\cal L}, \; s \in S_{\ell}, \; i \in {\cal I}_{\ell} \label{eq10}
\\
q_{i\ell s} \geq 0, \quad & \quad
\ell \in {\cal L}, \; s \in S_{\ell}, \; i \in {\cal I}_{\ell} \label{eq11a}
\\
y_{ij\ell s} \geq 0, \quad & \quad
\ell \in {\cal L}, \; s \in S_{\ell}, \; i,j \in {\cal I}_{\ell} \label{eq11b}
\\
I^+_{it}, \, I^-_{it} \geq 0, \quad & \quad
i \in {\cal I}, \; t \in {\cal T}. \label{eq11c}
\end{align}

Objective function~\eqref{fo} corresponds to the sum of the costs of inventory, backordering, setup, and production. Constraints~\eqref{eq2} define inventory and backordering conservation flow, i.e.\ the relation between inventory, backordering, demand, and produced quantities. Constraints~\eqref{eq3} determine that the total amount of inventory cannot exceed the warehousing capacity of the distribution center. In the problem under consideration, existing storage capacity cannot be increased nor violated. (If desired, it would be possible to consider, in the objective function, an additional cost associated with exceeding the current capacity.) Constraints~\eqref{eq4} guarantee that the sum of production time and setup time of each machine within each period is not larger than the corresponding machine's availability. Constraints~\eqref{eq5} ensure that machine~$\ell$ produces product~$i$ within subperiod~$s$ only if the machine has been previously configured for this purpose. Constraints~\eqref{eq6} impose a minimum lot size restriction. Constraints~\eqref{eq7} determine that, for every machine, within every machine's subperiod, a single product can be produced. Constraints~\eqref{eq8} indicate that, if the machine~$\ell$ switch from producing product~$i$ to producing product~$j$ occurs within subperiod~$s$, then the corresponding setup must be made over the course of subperiod~$s$. Constraints~(\ref{eq10},\ref{eq11a},\ref{eq11b},\ref{eq11c}) determine the variables' domain. Note that it is not necessary to impose variables~$y_{i j \ell s}$ to be binary, because they naturally assume values in $\{0,1\}$ at an optimal solution. This is due to the fact that constraints~\eqref{eq8} restrict them to be larger than or equal to $-1$, $0$, or~$1$; while constraints~\eqref{eq11b} inhibit the possibility of being negative. Then, the minimization of the objective function forces them to assume binary values at the optimal solution. 

\section{Relax-and-fix Heuristics} \label{secrelax}

\begin{figure}[!ht]
\centering
\includegraphics[scale=1.0]{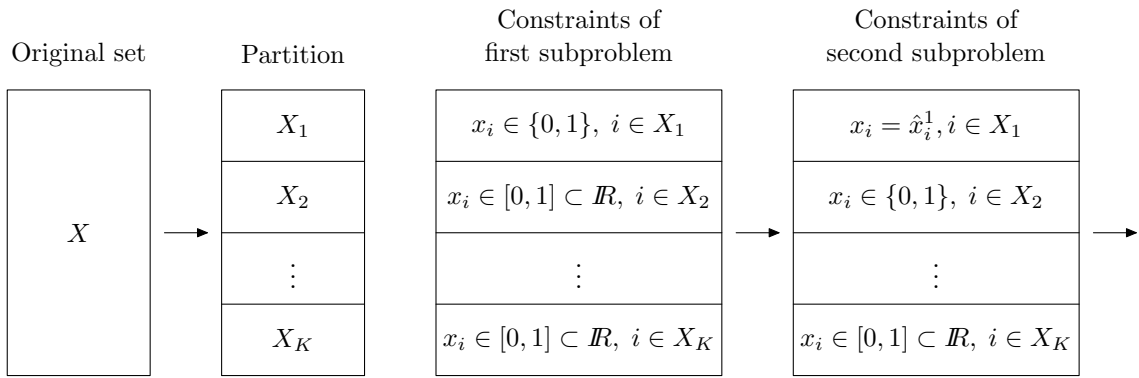}
\caption{Relax-and-fix algorithm working structure. In the second subproblem, $\hat x_i^1$, $i \in X_1$, correspond to the optimal values found when solving the first subproblem.}
\label{fig2}
\end{figure}

The main idea of the relax-and-fix (RF) heuristic \citep{pochetwolsey2006} for solving a MILP problem is to partition the set~$X$ of integer variables into $K$ subsets $X_1,\dots,X_K$ and to solve a sequence of $K$ smaller MILP problems. In the $k$th subproblem, the variables in $X_k$ are restricted to be integers; while the integrality of the variables in $X_{k+1} \cup \dots \cup X_K$ is relaxed and the variables in $X_1 \cup \dots \cup X_{k-1}$ are already fixed; see Figure~\ref{fig2}. 

Model~(\ref{fo}--\ref{eq11c}), that will be named ``Model ${\cal M}$'' from now on, needs to be modified in order to be used within the context of an RF-based heuristic. Let $X=\{ (i,\ell,s) \; | \; \ell \in {\cal L}, \; s \in S_{\ell}, \; i \in {\cal I}_{\ell} \}$ be the set of all valid indices' $3$-uples of variables $x_{i \ell s}$. The RF-based heuristic relies on the partition of the set $X$ into $K$ subsets $X_k$ ($k=1,\dots,K$) that verify $\cup_{k=1}^K X_k = X$ and $X_{k_1} \cap X_{k_2} = \emptyset$ for all $k_1 \neq k_2$. Let Model ${\cal M}_1$ be defined as Model~${\cal M}$ in which constraint~\eqref{eq10} is substituted with
\begin{equation}
\begin{array}{rl}
x_{i \ell s} \in \{0,1\}, & (i, \ell, s) \in X_1,\\[2mm]
0 \leq x_{i\ell s} \leq 1, & (i,\ell, s) \in \cup_{w=2}^K X_w.
\end{array}
\end{equation}
We can now define Model~${\cal M}_k$, for $k=2,\dots,K$, as Model~${\cal M}$ in which constraint~\eqref{eq10} is substituted with
\begin{align}
x_{i \ell s} = \hat x_{i \ell s}^w, & \quad (i, \ell, s) \in X_w \mbox{ for } w=1,\dots,k-1, \label{eq12} \\
x_{i \ell s} \in \{0,1\}, & \quad (i, \ell, s) \in X_k, \label{eq13} \\
0 \leq x_{i\ell s} \leq 1, & \quad (i,\ell, s) \in \cup_{w=k+1}^K X_w, \label{eq14}
\end{align}
where, in~\eqref{eq12}, $\hat x_{i \ell s}^w$ for $(i, \ell, s) \in X_w$ correspond to the optimal values obtained when solving Models~${\cal M}_w$ for $w=1,\dots,k-1$. This expresses the fact that, in Model ${\cal M}_k$, variables $x_{i \ell s}$ with $(i, \ell, s) \in \cup_{w=1}^{k-1} X_w$ are fixed, variables $x_{i \ell s}$ with $(i, \ell, s) \in X_k$ preserve their integrality constraint, and variables $x_{i \ell s}$ with $(i, \ell, s) \in \cup_{w=k+1}^{K} X_w$ are variables whose integrality constraint is relaxed. Note that, in Model~${\cal M}$, constraint~\eqref{eq10}, which says that $x_{i \ell s} \in \{0,1\}$ for all $(i, \ell, s) \in X$, implies that all variables~$y_{i j \ell s}$, for all $\ell \in {\cal L}$, $\in S_{\ell}$, $i, j \in {\cal I}_{\ell}$, assume binary values at an optimal solution. On the other hand, a different relation holds in Model~${\cal M}_k$ ($k=1,\dots,K$). Due to~\eqref{eq8}, a variable $y_{i j \ell s}$ is guaranteed to assume binary values in an optimal solution of Model~${\cal M}_k$ only if $(i, \ell, s-1)$ and $(j, \ell, s) \in X_k$.

The key feature of an RF-based heuristic is the determination of the number of subproblems~$K$ and the sets $X_1, \dots, X_K$. Problem-dependent and problem-independent strategies will be considered. The problem-dependent strategies divide variables $x_{i \ell s}$ considering the dimensions related to machine, product, and period; and they are based on several metrics associated with machines, products, and periods that we now describe. We define the demand~$d_i$ of a product~$i$ as
\begin{equation} \label{prodem}
d_i = \sum_{t \in {\cal T}} d_{it};
\end{equation}
and the flexibility $f_i$ of a product~$i$ as the number of machines that can produce it, i.e.\ 
\begin{equation} \label{proversa}
f_i = | {\cal L}_i |
\end{equation}
for all $i \in {\cal I}$. In addition, inspired by Vogel's approximation method for the transportation problem~\citep[\S4]{vogel} (see also~\cite{glover}), we define the discrepancy~$a_i$ of a product~$i$ as the difference between its two smallest processing times, given by
\begin{equation} \label{discrepancy}
a_i = p_{i \hat \ell_2} - p_{i \hat \ell_1}, \;
\hat \ell_1 = \argmin_{\ell \in {\cal L}_{i}} \{ p_{i\ell} \}, \;
\mbox{and }
\hat \ell_2 = \argmin_{\ell \in {\cal L}_{i} \setminus \{\hat \ell_1\}} \{ p_{i\ell} \}.
\end{equation}
Following \cite{dastidar2005scheduling}, we say a machine is critical if it is able, and thus it potentially will need to, process products with low flexibility. Therefore, we define the criticality~$\omega_\ell$ of a machine~$\ell$ as 
\begin{equation} \label{maqcrit}
\omega_\ell = m - \min_{i \in {\cal I}_\ell} \{ f_i \}
\end{equation}
for all $\ell \in {\cal L}$. We also associate with a machine~$\ell$ a metric of efficiency $\varepsilon_\ell$, given by the sum of the machine's average processing times and costs, i.e.\ 
\begin{equation} \label{maqeff}
\varepsilon_\ell = \frac{1}{| {\cal I}_\ell |} ( \sum_{{i \in \cal I}_\ell} p_{i\ell} + c_{i\ell}^P )
\end{equation}
for all $\ell \in {\cal L}$. Finally, for a period $t \in {\cal T}$, we define its overall demand $\delta_t$ as 
\begin{equation} \label{perdem}
\delta_t = \sum_{i \in {\cal I}} d_{it}.
\end{equation}

As a whole, nine different problem-dependent and two problem-independent partition strategies are considered. Each strategy consists of a way of sorting the 3-uples of indices $(i,\ell,s) \in X$, where $X$ is the set of valid 3-uples of the variables $x_{i \ell s}$. After sorting, the first $\lfloor |X|/K \rfloor$ variables constitute the set $X_1$, the next $\lfloor |X|/K \rfloor$ variables constitute the set $X_2$, and so on. In fact, if $|X|$ is not a multiple of $K$, then we have that $K \lfloor |X| / K \rfloor < |X|$. So, in the first $r = |X| - K \lfloor |X| / K \rfloor$ subsets, we consider $\lceil \cdot \rceil$ instead of $\lfloor \cdot \rfloor$ for the cardinality of $X_1, \dots, X_r$. In most strategies, the suggested order does not imply a total order of the 3-uples. Therefore, tie-breaking criteria play an important role in the strategies. In the suggested strategies, the tie-breaking criterion will always consist in applying a second strategy. Thus, the strategies are intrinsically hybrid. As a second tie-breaking rule, the lexicographic order of the index 3-uples is used. The description of each strategy follows:

\begin{itemize}
\item \textbf{Time-dimension-based strategies:}
    \begin{enumerate}[label=(S{{\arabic*}})]
    \item Chronological time: Given $(i, \ell, s) \in X$, we have that $s \in {\cal S}_{\ell t}$ for some $t \in {\cal T}$. 3-uples are sorted in increasing lexicographical order of $(t,s)$, where $t$ is the period~$t$ to which the subperiod $s$ corresponds.
    \item Periods with larger demand first: The same as (S1), but 3-uples $(i, \ell, s) \in X$ are sorted in decreasing order of~$\delta_t$ (instead of increasing order of~$t$), where~$\delta_t$ is the demand of period~$t$ given by~\eqref{perdem}.  
    \end{enumerate}

\item \textbf{The product-dimension-based strategies:}
    \begin{enumerate}[label=(S{{\arabic*}})]
    \setcounter{enumi}{2}
    \item Most demanded products first: 3-uples $(i, \ell, s) \in X$ are sorted in decreasing order of the items demands $d_i$ given by~\eqref{prodem}.
    \item Less demanded products first: The same as (S3) but in increasing order.
    \item Less flexible products first: 3-uples $(i, \ell, s) \in X$ are sorted in increasing order of the items flexibility $f_i$ given by~\eqref{proversa}.
    \item Products with more discrepant two smallest production times first: 3-uples $(i, \ell, s) \in X$ are sorted in decreasing order of the items discrepancy $a_i$ given by~\eqref{discrepancy}.
    \end{enumerate}
\item \textbf{The machine-dimension-based strategies:} 
    \begin{enumerate}[label=(S{{\arabic*}})]
    \setcounter{enumi}{6}
    \item Less efficient machines first: 3-uples $(i, \ell, s) \in X$ are sorted in increasing order of the machines efficiency $\varepsilon_{\ell}$ given by~\eqref{maqeff}.
    \item More efficient machines first: The same as (S7) but in decreasing order.
    \item More critical machines first: 3-uples $(i, \ell, s) \in X$ are sorted in decreasing order of the machines criticality $\omega_\ell$ given by~\eqref{maqcrit}.
    \end{enumerate}
\item \textbf{Problem-independent strategies:} 
    \begin{enumerate}[label=(S{{\arabic*}})]
    \setcounter{enumi}{9}
    \item More fractional variables first: In this strategy, we first solve a linear programming (LP) problem that corresponds to Model~${\cal M}$ in which the constraint $x_{i \ell s} \in \{0,1\}$ is relaxed to $0 \leq x_{i \ell s} \leq 1$ for \textit{all} $(i,\ell,s) \in X$. Let $\hat x_{i \ell s}^0$, for $(i,\ell,s) \in X$, be the optimal solution of the LP problem. For each variable $x_{i \ell s}$ for $(i,\ell,s) \in X$, we compute the ``distance to integrality'' given by
    \[
    d_{i \ell s}^0 = \min\{ \hat x_{i \ell s}^0, 1 - \hat x_{i \ell s}^0\}.
    \]
    3-uples $(i, \ell, s) \in X$ are sorted in decreasing order of their distance to integrality $d_{i \ell s}^0$. In fact, there is no need to sort all 3-uples since, to construct $X_1$, only the $\lfloor |X|/K \rfloor$ or $\lceil |X|/K \rceil$ 3-uples with largest $d_{i \ell s}^0$ are required. In general, after having solved the $k$th subproblem ($k<K$), distances $d_{i \ell s}^k= \min\{ \hat x_{i \ell s}^k, 1 - \hat x_{i \ell s}^k\}$ for $(i,\ell,s) \in X \setminus \cup_{w=1}^k X_w$ are computed, where $\hat x_{i \ell s}^k$ is the optimal solution of the $k$th subproblem, and the $\lfloor |X|/K \rfloor$ or $\lceil |X|/K \rceil$ 3-uples with largest $d_{i \ell s}^k$  are selected to constitute $X_{k+1}$.
    \item More influential variables first: In this strategy, variables with more ``influence'' in the objective function are considered first. This ``influence'' could be directly related to the cost associated with the variable in the objective function; and this is the most direct interpretation of this rule. However, in the problem at hand, variables $x_{i \ell s}$ with $(i,\ell,s) \in X$ do not appear in the objective function. However, their value influences the value of variables $q_{i \ell s}$ and $y_{ij\ell s}$. Thus, the influence $\alpha_{i \ell s}$ of variable $x_{i \ell s}$ is defined as
    \[
    \alpha_{i \ell s} = ( \sum_{j \in {\cal I}_{\ell}} c_{ij\ell}^S  ) + c_{i \ell}^P.
    \]
    Note that, in fact, $\alpha_{i \ell s}$ does not depend on $s$ and, therefore, variables $x_{i \ell s_1}$ and $x_{i \ell s_2}$ with $s_1 \neq s_2$ have the same influence measure. Therefore, in the particular problem at hand, this problem-independent strategy depends on the products and the machines. More specifically, products with large processing times and/or that, after being produced, require a time consuming preparation of the machine, are considered more influential.
    \end{enumerate}
\end{itemize}

Note that strategy (S10) is a dynamic strategy that differs of all other strategies because set $X_k$ is determined after having solved subproblems from~$1$ to $k-1$; while other strategies determine all subsets $X_1, \dots, X_K$ a priori. Moreover, strategy (S10) requires to solve an LP problem first. In the present work, we consider hybrid strategies that consists in applying a problem-dependent strategy between (S1) and (S9), using (S10) or (S11) as first tie-breaking rule, and the lexicographic order in the 3-uples $(i, \ell, s) \in X$ as second tie-breaking rule. (Hybrid strategies that combined two problem-dependent strategies were also evaluated numerically, but they showed marginal benefits in relation to those presented in this paper.)

\section{Numerical Experiments} \label{secnum}

In this section, we evaluate the performance of the partition strategies described in the previous section. In a first set of experiments, strategies are evaluated on a set of randomly generated instances. In a second experiment, selected strategies are applied to a set of ten real-world instances of the industry of personal care products. The experiments were carried out in a Intel Xeon X5690 3.47 GHz machine with 64 Gb of RAM. The RF algorithms and the formulations were solved using CPLEX 12.10.0 using default parameters with concert library and C++ programming language. The code was compiled using gcc 6.3.0 compiler with the Codeblocks 16.01 IDE. Benchmark instances and code are available at \url{https://github.com/kennedy94/GLSPPL-RF}.

\subsection{Experiments with randomly generated instances}

The benchmark test suite is composed of 25 randomly generated instances inspired by the production environment of the target industry. All instances have a planning horizon composed of 16 time periods divided into 7 subperiods each. As each period corresponds to one week, the planning horizon comprises approximately four months. (This does not mean that each subperiod corresponds to one day. This means that a machine can produce at most seven different products in a period of one week. The subperiods correspond to the production of different products on the machine and are of varying duration.) 

Five groups of instances were considered (A, B, C, D and E) and five different random instances were generated within each of the groups, with a total of 25 instances.  The random generation used either discrete or continuous uniform distributions, depending on the nature of each parameter. The number of machines $m$ for groups A, B, C, D, and E is 2, 3, 4, 5 and 7, respectively; while the number of products $n$ for groups A, B, C, D, and E is 8, 12, 16, 20 and 28, respectively. It is considered that the machines are not ready for the production of any item at the beginning of the planning horizon, i.e.\ $x_{i \ell 0}=0$ for all $i \in {\cal I}$ and $\ell \in {\cal L}_i$. Table~\ref{tab1} shows details of the random generation of all instances parameters. It should be noted that the random instances were generated from real-world data. In the table, $\bar e$ is the average of the setup times $e_{ij\ell}$ for $\ell \in {\cal L}$ and $i,j \in {\cal I}_\ell$; and $d_t$ (for $t \in {\cal T}$) corresponds to an aggregated period demand. For each item $i \in {\cal I}$, a proportion $d_{i}^p \in [0.05, 0.9]$ is randomly generated (independently of $t$); and $d_{it}$ is defined as $d_{it} = d_t \; d_i^p / ( \sum_{i \in {\cal I}} d_i^p )$. The production and setup time intervals for each group were obtained from the average times and standard deviations provided by the company. The range for the minimum lot $q_{i\ell}^{\mathrm{lb}}$ of product $i$ on machine $\ell$ relative to each group was built based on the average and the standard deviation, informed by the industrial area of the company, of the minimum feasible quantity to be produced with productivity and material loss within the required standards. In fact, average and standard deviation were informed for each group in working shifts $\hat q^{\mathrm{lb}}$ and transformed in units of products $q_{i\ell}^{\mathrm{lb}}$ considering the processing times~$p_{i\ell}$ and the fact that each working shift corresponds to 8 hours of labor. The ranges for the amount of products in stock and backordered products at the beginning of the time horizon were calculated using averages and standard deviations, based on historical data, provided by the company. Table~\ref{tab2} shows detailed information about the generated instances. In the table, columns IV, CV, and CO, correspond to the numbers of integer variables, continuous variables, and constraints of Model ${\cal M}$, respectively. The last line of each group displays average values. 

\begin{table}[ht!]
\centering
\resizebox{\textwidth}{!}{
\begin{tabular}{|c|c|c|c|c|c|c|}
\hline
\multirow{2}{*}{Parameter} & \multirow{2}{*}{Unit of measure} & \multicolumn{5}{c|}{Group} \\  \cline{3-7}
& & A & B & C & D & E \\
\hline
\hline
$m$ & unit & 2 & 3 & 4 & 5 & 7\\
$n$ & unit & 8 & 12 & 16 & 20 & 28\\
$|{\cal I}_{\ell}|$ & unit & $[5,8]$ & $[3,9]$ & $[4,10]$ & $[5,12]$ & $[2,12]$\\
$d_t$ & thousand units & $[9,13]$ & $[16,24]$ & $[19,78]$ & $[27,220]$ & $[65,100]$\\
$C^W$ & thousand units & $[10,14]$ & $[14,18]$ & $[40,48]$ & $[180,220]$ & $[120,150]$\\
$C^P_{\ell t}$ & hours & $160$ & $160$ & $160$ & $160$ & $160$\\
$h_i$ & cost per unit & $[0.27,0.54]$ & $[0.02,0.034]$ & $[0.03,0.08]$ & $[0.07,0.21]$ & $[0.087,0.433]$\\
$g_i$ & cost per unit & $[10,15] \times h_i$ & $[10,15] \times h_i$ & $[10,15] \times h_i$ & $[10,15] \times h_i$ & $[10,15] \times h_i$\\
$p_{i\ell}$ & hours & $[0.012,0.04]$ & $[0.008,0.05]$ & $[0.007,0.017]$ & $[0.003,0.01]$ & $[0.005,0.028]$\\
$c^P_{i \ell}$ & cost per unit & $[0.8,1.2] \times p_{i\ell}$ & $[0.8,1.2] \times p_{i\ell}$ & $[0.8,1.2] \times p_{i\ell}$ & $[0.8,1.2] \times p_{i\ell}$ & $[0.8,1.2] \times p_{i\ell}$\\
$\hat q^{lb}$ & working shifts & $[1,9]$ & $[3,6]$ & $[3,6]$ & $[6,9]$ & $[3,6]$\\
$e_{ij\ell}$ & hours & $[2,6]$ & $[5,9]$ & $[2,6]$ & $[1,6]$ & $[2,8]$\\
$c^S_{ij\ell}$ & cost per setup & $[80,100] \times e_{ij\ell}/\bar{e}$ & $[100,200] \times e_{ij\ell}/\bar{e}$ & $[100,200] \times e_{ij\ell}/\bar{e}$ & $[230,1200] \times e_{ij\ell}/\bar{e}$ & $[150,620] \times e_{ij\ell}/\bar{e}$\\
$I^+_{i0}$ & unit & $[0,4000]$ & $[0,4000]$ & $[0,20000]$ & $[0,50000]$ & $[0,30000]$\\
$I^-_{i0}$ & unit & $[0,500]$ & $[0,500]$ & $[0,2000]$ & $[0,4000]$ & $[0,4000]$\\
\hline
\end{tabular}}
\caption{Detailed information of the random generation of instances' parameters.}
\label{tab1}
\end{table}

\begin{table}[ht!]
\centering
\resizebox{\textwidth}{!}{
\begin{tabular}{|c|ccc|c|ccc|c|ccc|c|ccc|c|ccc|}
\hline
\multicolumn{4}{|c|}{Group A} &
\multicolumn{4}{c|}{Group B} &
\multicolumn{4}{c|}{Group C} &
\multicolumn{4}{c|}{Group D} &
\multicolumn{4}{c|}{Group E} \\
\hline
Inst. & IV & CV & CO & Inst. & IV & CV & CO  & Inst. & IV & CV & CO  & Inst. & IV & CV & CO  & Inst. & IV & CV & CO \\
\hline
\hline
A1    & 1,456 & 11,346 & 12,832 & B1    & 2,016 & 14,646 & 16,720 & C1    & 3,136 & 27,438 & 30,576 & D1    & 4,928 & 52,208 & 56,976 & E1    & 5,488 & 46,604 & 52,096 \\
A2    & 1,568 & 13,154 & 14,736 & B2    & 2,128 & 16,680 & 18,848 & C2    & 2,912 & 23,596 & 26,544 & D2    & 4,480 & 42,716 & 47,120 & E2    & 5,376 & 45,926 & 51,312 \\
A3    & 1,456 & 11,346 & 12,832 & B3    & 2,128 & 18,036 & 20,192 & C3    & 3,024 & 25,178 & 28,224 & D3    & 4,592 & 43,620 & 48,128 & E3    & 5,712 & 52,706 & 58,368 \\
A4    & 1,456 & 11,346 & 12,832 & B4    & 2,352 & 20,070 & 22,432 & C4    & 3,024 & 25,856 & 28,896 & D4    & 5,152 & 55,824 & 60,784 & E4    & 5,712 & 52,028 & 57,696 \\
A5    & 1,344 &  9,990 & 11,376 & B5    & 2,240 & 18,714 & 20,976 & C5    & 3,024 & 24,500 & 27,552 & D5    & 5,040 & 52,434 & 57,312 & E5    & 5,600 & 53,384 & 58,928 \\
\hline
avg.  & 1,456 & 11,436 & 12,922 & avg.  & 2,173 & 17,629 & 19,834 & avg.  & 3,024 & 25,314 & 28,358 & avg.  & 4,838 & 49,360 & 54,064 & avg.  & 5,578 & 50,130 & 55,680 \\
\hline                          
\end{tabular}}
\caption{Detailed information of benchmark random instances.}
\label{tab2}
\end{table}

Initially, to decide which of the problem-independent strategies would be used as the tie-breaking rule, we ran strategies S10 and S11 varying $K \in {1,2,\dots,10}$. The two rows at the bottom of Table~\ref{tab3} show the results. In the table, the results obtained with $K \geq 2$ are compared with the result obtained with $K=1$, which simply corresponds to solving Model ${\cal M}$ using CPLEX. The results appear separated by group and refer to the 5 instances of the group solved with 9 values of $K \in \{2,\dots,10\}$ . Under the heading $(W,G(\%))$ it appears how many, out of the 45 cases, the relax-and-fix strategy found a better result ($W$ stands for ``win'') than CPLEX and what the average gap was in these cases ($G(\%)$ stands for ``average gap in percentage''). Under the heading $(L,G(\%))$ appears how many of the $45$ cases did the relax-and-fix strategy perform worse ($L$ stands for ``lost'') than CPLEX and what was the average gap in these cases. Since there were no ties, if the number of wins and losses does not add up to 45, it is because for some instances and values of $K$, relax-and-fix could not find a feasible solution. The last column of the table shows the average gap considering the $25$ instances and the 9 tested values of $K$. The results show that the strategy S11 obtained better results than the strategy S10 and, for this reason, it will be used as a tie-breaker for all the other strategies that depends on the problem. The performance of strategies S1 to S9, using S11 as the tiebreaker strategy is shown at the top of Table~3. At this point it is important to mention that for both CPLEX and the relax-and-fix strategy a time limit of~1 hour was used. (The influence of this time limit on the comparison will be analyzed later in this section.) Furthermore, in the relax-and-fix strategy, the time was divided linearly between the subproblems such that the first subproblem has twice as much time as the last. Preliminary tests giving three times as much time to the first subproblem and dividing the time evenly among the subproblems showed similar results. The last column of the table shows that, in aggregate, that is, without evaluating the different values of~$K$ and the groups of instances individually, all strategies improved the results obtained with CPLEX with a gap ranging from $-36\%$ to $-47\%$. The numbers also show that strategies S1 and S9 presented the best results. It is also worth noting that these two strategies found feasible solutions for all instances and all tested values of~$K$. 

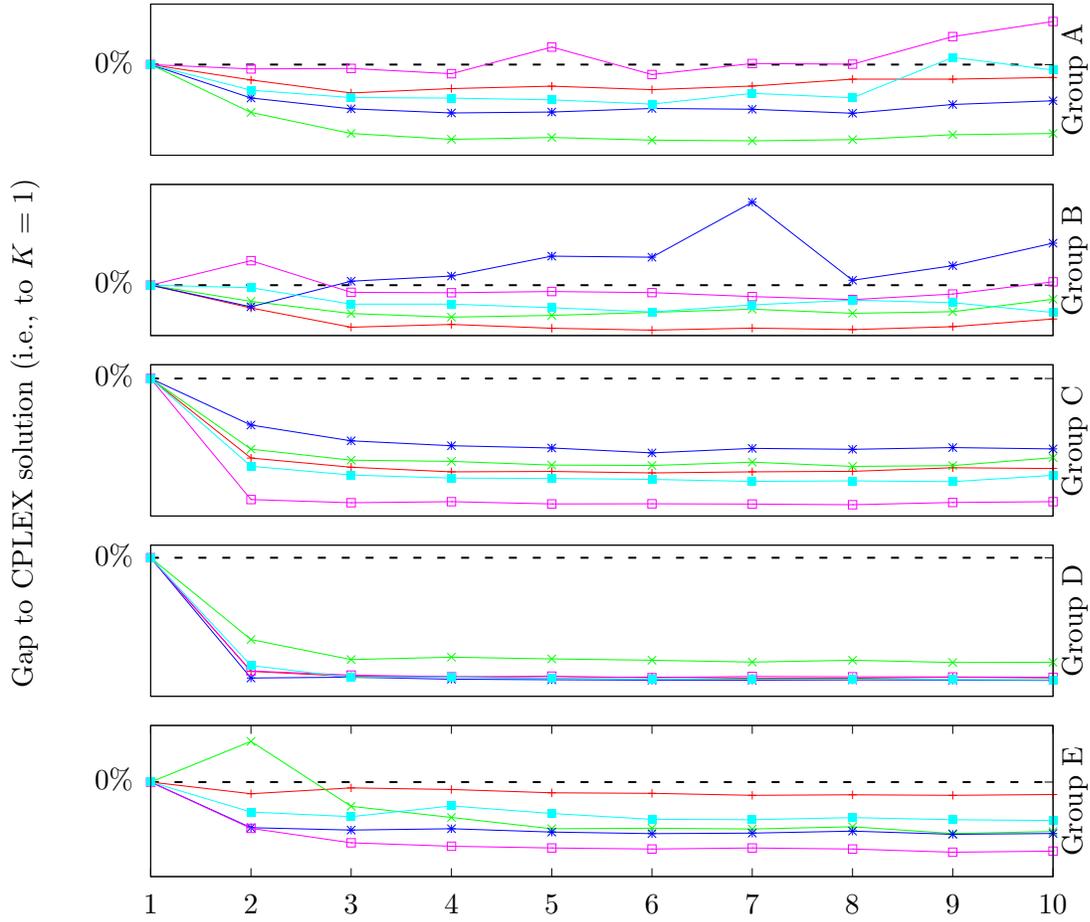
\begin{figure}[!ht]
\centering
\begin{center}
\input{abkm2020fig3.tex}
\end{center}
\caption{Gaps to CPLEX solution (i.e.\ $K=1$) of strategy S1 on random benchmark instances.}
\label{fig3}
\end{figure}

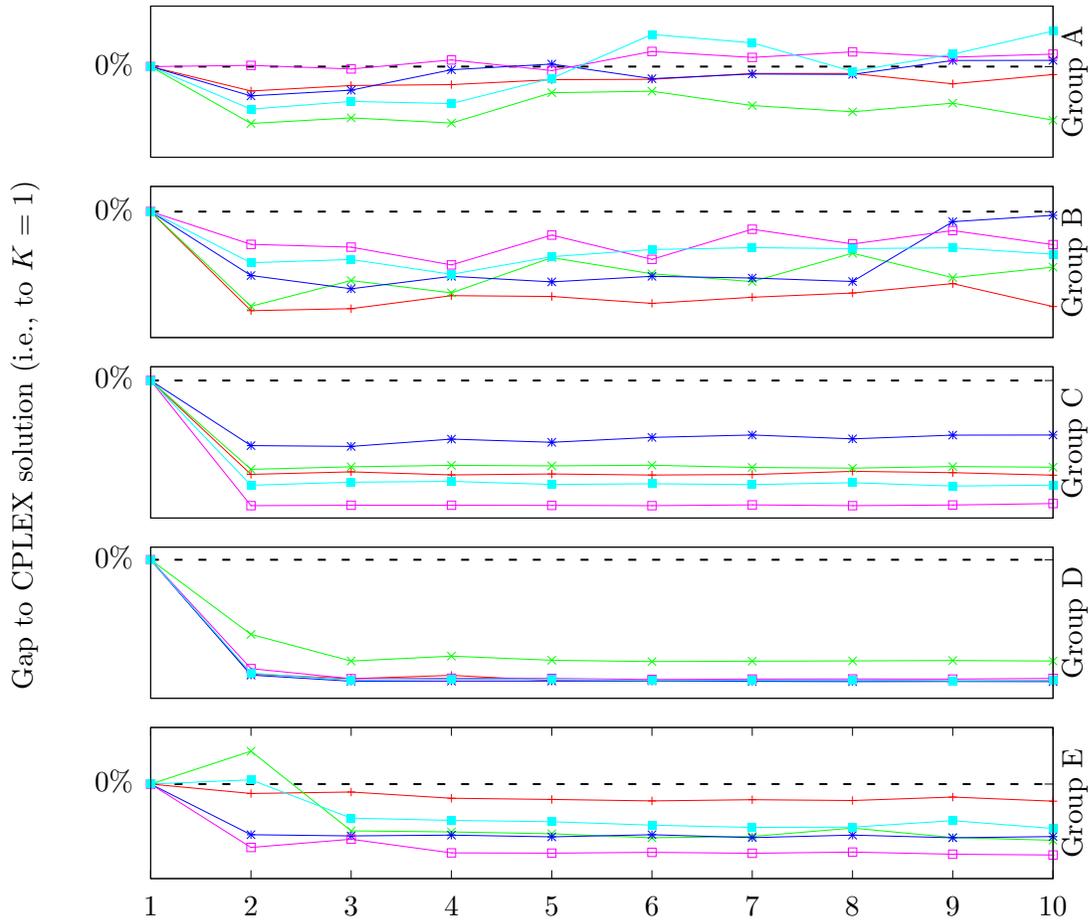
\begin{figure}[!ht]
\centering
\begin{center}
\input{abkm2020fig4.tex}
\end{center}
\caption{Gaps to CPLEX solution (i.e. $K=1$) of strategy S9 on random benchmark instances.}
\label{fig4}
\end{figure}

Table~\ref{tab4} shows the results for strategies S1 and S9 disaggregated, that is, for each value of~$K$ separately. Figures~\ref{fig3} and~\ref{fig4} show graphically the same that is being shown by Table~\ref{tab4}. The numbers in the table show that, except for the case $K=2$, the results vary little depending on the value of $K$ chosen, which can be seen as a positive feature of the methods. The numbers also show that the relax-and-fix strategies almost always beat CPLEX, except for some instances of groups~A and~B that concentrate the smallest instances.

\begin{table}[!ht]
\begin{center}
\resizebox{\textwidth}{!}{
\begin{tabular}{|c|rrrr|rrrr|rrrr|rrrr|rrrr|r|}
\hline
\multicolumn{1}{|c|}{\multirow{2}{*}{Strategy}} &
\multicolumn{4}{|c}{Group A} & 
\multicolumn{4}{|c}{Group B} & 
\multicolumn{4}{|c}{Group C} & 
\multicolumn{4}{|c}{Group D} & 
\multicolumn{4}{|c}{Group E} & 
\multicolumn{1}{|c|}{\multirow{2}{*}{G(\%)}}\\
\cline{2-21}
& \multicolumn{4}{|c}{(W,G(\%)) / (L,G(\%))} 
& \multicolumn{4}{|c}{(W,G(\%)) / (L,G(\%))} 
& \multicolumn{4}{|c}{(W,G(\%)) / (L,G(\%))} 
& \multicolumn{4}{|c}{(W,G(\%)) / (L,G(\%))} 
& \multicolumn{4}{|c}{(W,G(\%)) / (L,G(\%))} 
& \multicolumn{1}{|c|}{}\\
\hline
\hline
S1 & 39 & -24.71 & 6 & 10.71 & 35 & -23.69 & 10 & 24.57 & 45 & -67.39 & \multicolumn{2}{c|}{--} & 45 & -91.92 & \multicolumn{2}{c|}{--} & 44 & -42.56 & 1 & 43.22 & -46.58\\
S2 & 34 & -19.55 & 11 & 11.77 & 35 & -20.76 & 10 & 19.60 & 45 & -62.98 &  \multicolumn{2}{c|}{--} & 44 & -91.14 & \multicolumn{2}{c|}{--} & 43 & -39.35 & 2 & 15.36 & -42.73\\
S3 & 33 & -15.09 & 12 & 12.14 &  37 & -19.43 & 8 & 11.94 & 45 & -62.48 & \multicolumn{2}{c|}{--} &  44  & -91.13 & \multicolumn{2}{c|}{--} &  40 & -42.19 & 5 & 21.80 & -41.85
\\
S4 & 34 & -20.92 & 11 & 5.61 & 40 & -20.73 & 5 & 14.17 & 45 & -62.91 & \multicolumn{2}{c|}{--} & 45 & -90.34 & \multicolumn{2}{c|}{--} & 39& -40.46 & 6 & 39.41 & -42.87\\
S5 & 28 & -16.21 & 17 & 14.26 & 32 & -22.93 & 13 & 7.55 & 45 & -67.02 & \multicolumn{2}{c|}{--} & 45 & -90.54 & \multicolumn{2}{c|}{--} & 44 & -43.31 & 1 & 38.10 & -43.58
\\
S6 &
35 & -20.60 & 10 & 6.46 &
37 & -22.34 & 8 & 8.64 &
45 & -63.53 & \multicolumn{2}{c|}{--} &
43 & -89.67 & 2 & 3.87 &
43 & -41.15 & 2 & 49.82 &
-43.95\\
S7 &
24 & -12.53 & 21 & 10.47 & 32 & -19.92 & 13 & 37.50 & 44 & -60.14 & 1 & 21.70 & 45 & -88.13 & \multicolumn{2}{c|}{--} & 40 & -36.75 & 5 & 22.77 & -36.34\\
S8 &
27 & -18.42 & 18 & 8.77 & 23 & -20.01 & 22 & 15.84 & 45 & -64.13 & \multicolumn{2}{c|}{--} & 45 & -89.75 & \multicolumn{2}{c|}{--} & 43 & -41.05 & 2 & 18.19 & -40.46
\\
S9 &
31 & -15.71 & 14 & 8.85 & 45 & -22.58 & \multicolumn{2}{c|}{--} & 45 & -67.89 & \multicolumn{2}{c|}{--} & 45 & -91.38 & \multicolumn{2}{c|}{--} & 43 &-47.13 & 2 & 39.43 & -46.81\\
\hline
S10 &
23 & -15.73 & 22 & 25.73 &
16 & -12.79 & 29 & 20.25 &
45 & -59.67 & \multicolumn{2}{c|}{--} &
42 & -89.89 & 1 & 1.09 &
42 & -31.93 & 3 & 40.28 &
-31.93
\\
S11 &
31 & -17.01 & 14 & 9.04 &
45 & -21.34 & \multicolumn{2}{c|}{--} &
45 & -66.73 & \multicolumn{2}{c|}{--} &
45 & -91.84 & \multicolumn{2}{c|}{--} &
43 & -38.79 & 2 & 56.22 &
-44.92
\\
\hline
\end{tabular}}
\end{center}
\caption{General performance of hybrid strategies S1 to S9 with S11 as a tie-breaking rule plus standalone strategies S10 and S11 on random instances.}
\label{tab3}
\end{table}

\begin{table}[!ht]
\begin{center}
\resizebox{\textwidth}{!}{
\begin{tabular}{|c|rrrr|rrrr|rrrr|rrrr|rrrr|r|}
\hline
\multicolumn{22}{|c|}{Strategy S1}\\
\hline
\multicolumn{1}{|c|}{\multirow{2}{*}{$K$}} &
\multicolumn{4}{|c}{Group A} & 
\multicolumn{4}{|c}{Group B} & 
\multicolumn{4}{|c}{Group C} & 
\multicolumn{4}{|c}{Group D} & 
\multicolumn{4}{|c}{Group E} & 
\multicolumn{1}{|c|}{\multirow{2}{*}{G(\%)}}\\
\cline{2-21}
& \multicolumn{4}{|c}{(W,G(\%)) / (L,G(\%))} 
& \multicolumn{4}{|c}{(W,G(\%)) / (L,G(\%))} 
& \multicolumn{4}{|c}{(W,G(\%)) / (L,G(\%))} 
& \multicolumn{4}{|c}{(W,G(\%)) / (L,G(\%))} 
& \multicolumn{4}{|c}{(W,G(\%)) / (L,G(\%))} 
& \multicolumn{1}{|c|}{}\\
\hline
\hline
2 & 5 & -16.77 & \multicolumn{2}{c|}{--} & 4 & -15.67 & 1 & \multicolumn{1}{c|}{24.35} & 5 & -58.99 & \multicolumn{2}{c|}{--} & 5 & -85.23 & \multicolumn{2}{c|}{--} & 4 & -35.27 & 1 & \multicolumn{1}{c|}{43.20} & -37.65
 \\
3 & 5 & -23.61 & \multicolumn{2}{c|}{--} & 4 & -23.93 & 1 & \multicolumn{1}{c|}{4.00} & 5 & -65.99 & \multicolumn{2}{c|}{--} & 5 & -91.60 & \multicolumn{2}{c|}{--} & 5 & -36.63 & \multicolumn{2}{c|}{--} & -47.24 \\
4 & 5 & -25.13 & \multicolumn{2}{c|}{--} & 4 & -24.25 & 1 & \multicolumn{1}{c|}{9.10} & 5 & -67.81 & \multicolumn{2}{c|}{--} & 5 & -91.85 & \multicolumn{2}{c|}{--} & 5 & -37.53 & \multicolumn{2}{c|}{--} & -47.98 \\
5 & 4 & -29.32 & 1 & \multicolumn{1}{c|}{11.54} & 4 & -25.37 & 1 & \multicolumn{1}{c|}{28.81} & 5 & -69.01 & \multicolumn{2}{c|}{--} & 5 & -92.48 & \multicolumn{2}{c|}{--} & 5 & -43.25 & \multicolumn{2}{c|}{--} & -48.09 \\
6 & 5 & -25.67 & \multicolumn{2}{c|}{--} & 4 & -26.36 & 1 & \multicolumn{1}{c|}{27.80} & 5 & -70.11 & \multicolumn{2}{c|}{--} & 5 & -93.02 & \multicolumn{2}{c|}{--} & 5 & -45.05 & \multicolumn{2}{c|}{--} & -49.88  \\
7 & 4 & -28.31 & 1 & \multicolumn{1}{c|}{0.66} & 4 & -24.35 & 1 & \multicolumn{1}{c|}{82.37} & 5 & -69.15 & \multicolumn{2}{c|}{--} & 5 & -93.25 & \multicolumn{2}{c|}{--} & 5 & -45.41 & \multicolumn{2}{c|}{--} & -46.67  \\
8 & 4 & -28.35 & 1 & \multicolumn{1}{c|}{0.29} & 4 & -25.29 & 1 & \multicolumn{1}{c|}{4.93} & 5 & -69.83 & \multicolumn{2}{c|}{--} & 5 & -93.02 & \multicolumn{2}{c|}{--} & 5 & -44.17 & \multicolumn{2}{c|}{--} & -49.78  \\
9 & 3 & -27.46 & 2 & \multicolumn{1}{c|}{11.65} & 4 & -23.43 & 1 & \multicolumn{1}{c|}{19.30} & 5 & -68.73 & \multicolumn{2}{c|}{--} & 5 & -93.32 & \multicolumn{2}{c|}{--} & 5 & -47.52 & \multicolumn{2}{c|}{--} & -47.25 \\
10 & 4 & -20.42 & 1 & \multicolumn{1}{c|}{28.48} & 3 & -24.87 & 2 & \multicolumn{1}{c|}{22.52} & 5 & -66.87 & \multicolumn{2}{c|}{--} & 5 & -93.53 & \multicolumn{2}{c|}{--} & 5 & -46.71 & \multicolumn{2}{c|}{--} & -44.73  \\
\hline
\multicolumn{22}{c}{}\\
\hline
\multicolumn{22}{|c|}{Strategy S9}\\
\hline
\multicolumn{1}{|c|}{\multirow{2}{*}{$K$}} &
\multicolumn{4}{|c}{Group A} & 
\multicolumn{4}{|c}{Group B} & 
\multicolumn{4}{|c}{Group C} & 
\multicolumn{4}{|c}{Group D} & 
\multicolumn{4}{|c}{Group E} & 
\multicolumn{1}{|c|}{\multirow{2}{*}{G(\%)}}\\
\cline{2-21}
& \multicolumn{4}{|c}{(W,G(\%)) / (L,G(\%))} 
& \multicolumn{4}{|c}{(W,G(\%)) / (L,G(\%))} 
& \multicolumn{4}{|c}{(W,G(\%)) / (L,G(\%))} 
& \multicolumn{4}{|c}{(W,G(\%)) / (L,G(\%))} 
& \multicolumn{4}{|c}{(W,G(\%)) / (L,G(\%))} 
& \multicolumn{1}{|c|}{}\\
\hline
\hline
2 & 4 & -25.31 & 1 & 0.76 & 5 & -27.11 & \multicolumn{2}{c|}{--} & 5 & -69.48 & \multicolumn{2}{c|}{--} & 5 & -83.61 & \multicolumn{2}{c|}{--} & 3 & -43.51 & 2 & 19.72 & {-43.70} \\
3 & 5 & -17.37 & \multicolumn{2}{c|}{--} & 5 & -25.93 & \multicolumn{2}{c|}{--} & 5 & -68.37 & \multicolumn{2}{c|}{--} & 5 & -92.26 & \multicolumn{2}{c|}{--} & 5 & -41.38 & \multicolumn{2}{c|}{--} & {-49.06} \\
4 & 4 & -18.93 & 1 & {4.41} & 5 & -27.49 & \multicolumn{2}{c|}{--} & 5 & -67.42 & \multicolumn{2}{c|}{--} & 5 & -90.69 & \multicolumn{2}{c|}{--} & 5 & -46.15 & \multicolumn{2}{c|}{--} & {-49.20} \\
5 & 4 & -9.06 & 1 & {1.64} & 5 & -21.42 & \multicolumn{2}{c|}{--} & 5 & -68.32 & \multicolumn{2}{c|}{--} & 5 & -92.29 & \multicolumn{2}{c|}{--} & 5 & -47.49 & \multicolumn{2}{c|}{--} & {-47.29} \\
6 & 3 & -10.89 & 2 & {15.60} & 5 & -24.17 & \multicolumn{2}{c|}{--} & 5 & -67.59 & \multicolumn{2}{c|}{--} & 5 & -92.82 & \multicolumn{2}{c|}{--} & 5 & -48.72 & \multicolumn{2}{c|}{--} & {-46.72} \\
7 & 3 & -11.75 & 2 & {10.89} & 5 & -21.90 & \multicolumn{2}{c|}{--} & 5 & -67.50 & \multicolumn{2}{c|}{--} & 5 & -92.67 & \multicolumn{2}{c|}{--} & 5 & -49.50 & \multicolumn{2}{c|}{--} & {-46.85} \\
8 & 4 & -10.69 & 1 & {9.70} & 5 & -20.82 & \multicolumn{2}{c|}{--} & 5 & -67.54 & \multicolumn{2}{c|}{--} & 5 & -92.58 & \multicolumn{2}{c|}{--} & 5 & -47.04 & \multicolumn{2}{c|}{--} & {-46.92} \\
9 & 2 & -17.82 & 3 & {6.16} & 5 & -16.15 & \multicolumn{2}{c|}{--} & 5 & -67.33 & \multicolumn{2}{c|}{--} & 5 & -92.76 & \multicolumn{2}{c|}{--} & 5 & -47.95 & \multicolumn{2}{c|}{--} & {-45.52} \\
10 & 2 & -20.33 & 3 & {11.97} & 5 & -18.24 & \multicolumn{2}{c|}{--} & 5 & -67.44 & \multicolumn{2}{c|}{--} & 5 & -92.68 & \multicolumn{2}{c|}{--} & 5 & -50.94 & \multicolumn{2}{c|}{--} & -46.05 \\
\hline
\end{tabular}}
\end{center}
\caption{Detailed performance of strategies S1 and S9 with S11 as a tie-breaking rule varying $K \in \{2, 3, \dots, 10\}$ on random instances.}
\label{tab4}
\end{table}

The comparison with CPLEX depends very much on the 1 hour time limit being used; since if the instance is small and CPLEX is able to find an optimal solution (regardless of whether it can prove that the solution is optimal or not), then there is nothing that relax-and-fix can do. So we decided to test the influence of the time limit on the comparison. To this end, we re-run CPLEX and strategies S1 and S9 with time limits of $600$, $1200$, $1800$, $2400$, $3000$ and $3600$ seconds. We considered S1 with $K=6$ and S9 with $K=4$ because Table~\ref{tab4} shows that these were the best values of~$K$ for these two strategies. But this does not mean that we intend to use these values of~$K$ in the next experiments, since as already mentioned before, the method is robust and shows small variations with respect to the value of~$K$. Figures~\ref{fig5} and~\ref{fig6} show the results for strategies S1 (with $K=6$) and S9 (with $K=4$), respectively. Figures~\ref{fig5}(a) and~\ref{fig6}(a) show that when the time limit is reduced the advantage of the relax-and-fix strategy increases. In the case of strategy S1, the average gap with CPLEX, which is $-49.88\%$ for the one-hour time limit, goes to $-76.00\%$.  In the case of strategy S9, the gap goes from $-49.20\%$ to $-70.21\%$. This experiment shows that if you have less time, using the relax-and-fix strategy is even more advantageous. But the question that remains is: With less time, how much do the solutions found by relax-and-fix deteriorate? Figures~\ref{fig5}(b) and~\ref{fig6}(b) compare the solution obtained by relax-and-fix with a reduced budget against the solution found with the one hour budget. The figures show that the maximum deterioration is, in average, approximately $30\%$. The deterioration of the solutions found by CPLEX with reduced time is much greater, which is why the advantage of the relax-and-fix strategy increases.

\begin{figure}[!ht]
\centering
\begin{tabular}[!ht]{c}
\input{abkm2020fig5a.tex}\\
(a)\\
\input{abkm2020fig5b.tex}\\
(b)
\end{tabular}
\caption{(a) Boxplot of the gaps between CPLEX and strategy S1 ($K=6$) applied to the~25 random instances varying the time limits of the methods. (The average gaps are $-76.00\%$, $-65.45\%$, $-57.14\%$, $-50.82\%$, $-52.01\%$, $-49.88\%$ for the time limits $600, 1200, \dots 3600$, respectively). (b) Boxplot of the gaps between S1 ($K=6$) applied to the 25 random instances varying its time limit in $\{600, 1200, \dots, 3000\}$ and S1 ($K=6$) with a one-hour time limit. (The average gaps are $30.91\%$, $16.22\%$, $8.27\%$, $5.33\%$, and $2.43\%$ for the time limits $600, 1200, \dots 3000$, respectively.)}
\label{fig5}
\end{figure}

\begin{figure}[!ht]
\centering
\begin{tabular}[!ht]{c}
\input{abkm2020fig6a.tex}\\
(a)\\
\input{abkm2020fig6b.tex}\\
(b)
\end{tabular}
\caption{(a) Boxplot of the gaps between CPLEX and strategy S9 ($K=4$) applied to the~25 random instances varying the time limits of the methods. (The average gaps are $-70.21\%$, $-58.54\%$, $-53.69\%$, $-48.84\%$, $-50.73\%$, $-49.20\%$ for the time limits $600, 1200, \dots 3600$, respectively). (b) Boxplot of the gaps between S9 ($K=4$) applied to the~$25$ random instances varying its time limit in $\{600, 1200, \dots, 3000\}$ and S9 ($K=4$) with a one-hour time limit. (The average gaps are $34.18\%$, $87.55\%$, $12.87\%$, $9.34\%$, and $5.98\%$ for the time limits $600, 1200, \dots 3000$, respectively.)}
\label{fig6}
\end{figure}

\subsection{Experiments with real-world instances}

In this section we apply the presented methods to eight real instances provided by the company. All instances, as well as the randomly generated ones, have~16 periods divided into~7 sub-periods. The number of machines varies between~2 and~7 and the number of products varies between~8 and~26. Table~\ref{tab5} shows some details of the instances and their respective models. Recall that in the table, IV, CV and CO stands for ``integer variables'', ``continuous variables'', and ``constraints''. The table also shows, for each instance, the solution found by the company, which was calculated by an expert with an undisclosed empirical method.

\begin{table}[ht!]
\centering
\begin{tabular}{|c|cc|ccc|r|}
\hline
Inst. & $m$ & $n$ & IV & CV & CO & Company's solution\\
\hline
\hline
P1 & 4 &  9 & 2,016 & 12,058 & 14,336 & 1,069,419\\
P2 & 3 & 12 & 2,016 & 14,872 & 16,944 &    64,706\\
P3 & 4 &  8 & 2,128 & 14,284 & 16,672 &   754,967\\
P4 & 5 & 13 & 2,800 & 20,556 & 23,600 &   888,172\\
P5 & 2 & 20 & 3,584 & 64,186 & 67,120 &    51,740\\
P6 & 5 & 24 & 4,480 & 45,338 & 49,648 &   903,501\\
P7 & 7 & 26 & 5,040 & 40,660 & 45,792 &   636,216\\
P8 & 7 & 26 & 5,264 & 56,028 & 61,248 & 2,301,544\\
\hline
\end{tabular}
\caption{Detailed information of the real-world instances.}
\label{tab5}
\end{table}

\begin{figure}
\centering
\begin{tabular}{cccc}
\hspace{1cm} & \input{abkm2020fig7a.tex} & & \input{abkm2020fig7b.tex}\\
\\
\\
& (a) S1 & & (b) S9
\end{tabular}
\caption{Average gap to the company's solution of solutions found by strategies S1 and S9 varying $K \in \{2, 3, \dots, 10 \}$ and the CPU time limit in $\{600, 1200, \dots, 3600\}$ seconds.}
\label{fig7}
\end{figure}
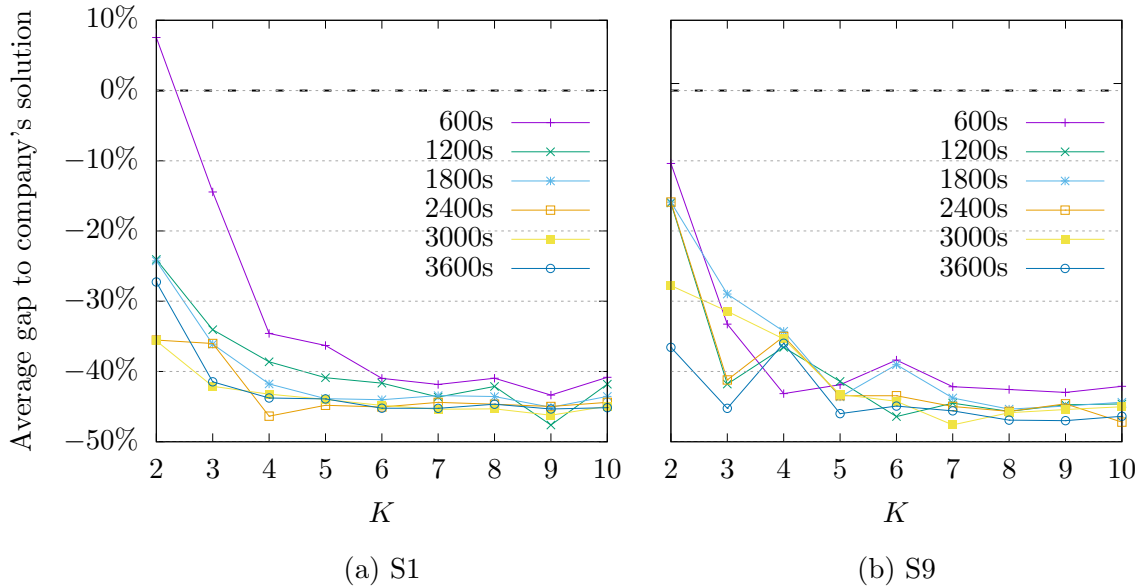

The results with the random instances showed that there is no clear advantage of a certain value of $K$ over others (excluding small values of $K$) and that small variations in the values of~$K$ generate small variations in the results. Because of this and also because we do not know how the solutions given by the company were calculated, we solved the eight instances with strategies~S1 and~S9 varying the time limit and varying~$K$, as we did with the random instances. Figure~\ref{fig7} shows the results. On the one hand, the figure shows that, for fixed~$K$, the longer the time limit, the better the solution found. On the other hand, the figure also shows that the two strategies find very similar values for any~$K \geq 6$. Because of this, we arbitrarily fixed, for both strategies, $K=8$.

Table~\ref{tab6} shows the results obtained by strategies~S1 and~S9, both with $K=8$, varying the time limit. The table also shows the results obtained by CPLEX, varying the time limit as well. For each method and time limit, the table shows the solution obtained and also shows the gap in relation to the solution presented by the company, calculated as
\[
100 \times \left( ( F_{\mathrm{method}} - F_{\mathrm{company}} ) / F_{\mathrm{company}} \right) \%,
\]
where $F_{\mathrm{method}}$ and $F_{\mathrm{company}}$ correspond to the objective function values of the solutions found by the method and the company, respectively. The numbers in the table show that the two proposed strategies always improve the company's solution by an amount that roughly varies, on average, between 41\% and 47\%, depending on the time limit given to the strategy. The solutions found by CPLEX are, on average, worse than the solutions reported by the company when the timeout is 600 or 1200 seconds; while CPLEX improves the company's solutions on average between 24\% and 35\% for longer times limits.

\begin{table}[ht!]
\begin{center}
\resizebox{\textwidth}{!}{%
\begin{tabular}{|c|rr|rr|rr|rr|rr|rr|}
\hline
\multicolumn{13}{|c|}{CPLEX}\\
\hline
\multirow{3}{*}{Inst.} &  
\multicolumn{2}{c|}{600s} &  
\multicolumn{2}{c|}{1200s} & 
\multicolumn{2}{c|}{1800s} & 
\multicolumn{2}{c|}{2400s} & 
\multicolumn{2}{c|}{3000s} & 
\multicolumn{2}{c|}{3600s}\\ 
\cline{2-13}
& \multicolumn{1}{c}{Objective} & \multicolumn{1}{c|}{\multirow{2}{*}{gap (\%)}} 
& \multicolumn{1}{c}{Objective} & \multicolumn{1}{c|}{\multirow{2}{*}{gap (\%)}} 
& \multicolumn{1}{c}{Objective} & \multicolumn{1}{c|}{\multirow{2}{*}{gap (\%)}} 
& \multicolumn{1}{c}{Objective} & \multicolumn{1}{c|}{\multirow{2}{*}{gap (\%)}} 
& \multicolumn{1}{c}{Objective} & \multicolumn{1}{c|}{\multirow{2}{*}{gap (\%)}} 
& \multicolumn{1}{c}{Objective} & \multicolumn{1}{c|}{\multirow{2}{*}{gap (\%)}} \\
& \multicolumn{1}{c}{function value} &
& \multicolumn{1}{c}{function value} &
& \multicolumn{1}{c}{function value} &
& \multicolumn{1}{c}{function value} &
& \multicolumn{1}{c}{function value} &
& \multicolumn{1}{c}{function value} &\\
\hline
\hline
P1 & 643,694.00 & -39.81 & 634,612.00 & -40.66 & 634,373.00 & -40.68 & 636,119.00 & -40.52 & 637,932.00 & -40.35 & 635,916.00 & -40.54\\
P2  & 22,932.20 & -64.56 & 22,495.10 & -65.23 & 20,988.70 & -67.56 & 20,690.20 & -68.02 & 20,775.50 & -67.89 & 20,753.30 & -67.93\\
P3 & 479,236.00 & -36.52 & 475,105.00 & -37.07 & 474,926.00 & -37.09 & 474,469.00 & -37.15 & 474,895.00 & -37.10 & 473,315.00 & -37.31\\
P4 & 485,099.00 & -45.38 & 482,778.00 & -45.64 & 480,030.00 & -45.95 & 460,206.00 & -48.19 & 457,134.00 & -48.53 & 450,809.00 & -49.24\\
P5 & 42,116.70 & -18.60 & 34,643.20 & -33.04 & 34,643.20 & -33.04 & 34,643.20 & -33.04 & 34,643.20 & -33.04 & 34,643.20 & -33.04\\
P6 & 2,454,100.00 & 171.62 & 1,073,460.00 & 18.81 & 832,150.00 & -7.90 & 832,150.00 & -7.90 & 832,150.00 & -7.90 & 870,273.00 & -3.68\\
P7 & 498,090.00 & -21.71 & 413,557.00 & -35.00 & 422,536.00 & -33.59 & 417,686.00 & -34.35 & 414,406.00 & -34.86 & 416,453.00 & -34.54\\
P8 & 14,885,400.00 & 546.76 & 14,885,400.00 & 546.76 & 3,935,990.00 & 71.02 & 1,975,860.00 & -14.15 & 1,975,860.00 & -14.15 & 1,975,860.00 & -14.15\\ 
\hline
&  & 61.48 &  & 38.62 &  & -24.35 &  & -35.42 &  & -35.48 &  & -35.05\\
\hline
\multicolumn{13}{c}{}\\
\hline
\multicolumn{13}{|c|}{Strategy S1 ($K = 8$)}\\
\hline
\multirow{3}{*}{Inst.} &  
\multicolumn{2}{c|}{600s} &  
\multicolumn{2}{c|}{1200s} & 
\multicolumn{2}{c|}{1800s} & 
\multicolumn{2}{c|}{2400s} & 
\multicolumn{2}{c|}{3000s} & 
\multicolumn{2}{c|}{3600s}\\ 
\cline{2-13}
& \multicolumn{1}{c}{Objective} & \multicolumn{1}{c|}{\multirow{2}{*}{gap (\%)}} 
& \multicolumn{1}{c}{Objective} & \multicolumn{1}{c|}{\multirow{2}{*}{gap (\%)}} 
& \multicolumn{1}{c}{Objective} & \multicolumn{1}{c|}{\multirow{2}{*}{gap (\%)}} 
& \multicolumn{1}{c}{Objective} & \multicolumn{1}{c|}{\multirow{2}{*}{gap (\%)}} 
& \multicolumn{1}{c}{Objective} & \multicolumn{1}{c|}{\multirow{2}{*}{gap (\%)}} 
& \multicolumn{1}{c}{Objective} & \multicolumn{1}{c|}{\multirow{2}{*}{gap (\%)}} \\
& \multicolumn{1}{c}{function value} &
& \multicolumn{1}{c}{function value} &
& \multicolumn{1}{c}{function value} &
& \multicolumn{1}{c}{function value} &
& \multicolumn{1}{c}{function value} &
& \multicolumn{1}{c}{function value} &\\
\hline
\hline
P1 & 652,170.00 & -39.02 & 648,090.00 & -39.40 & 647,849.00 & -39.42 & 647,849.00 & -39.42 & 644,788.00 & -39.71 & 642,909.00 & -39.88\\
P2  & 23,447.80 & -63.76 & 23,575.70 & -63.56 & 23,536.40 & -63.63 & 23,588.40 & -63.55 & 23,575.70 & -63.56 & 20,670.90 & -68.05\\
P3 & 481,900.00 & -36.17 & 478,173.00 & -36.66 & 476,597.00 & -36.87 & 476,597.00 & -36.87 & 476,764.00 & -36.85 & 472,957.00 & -37.35\\
P4 & 467,198.00 & -47.40 & 456,614.00 & -48.59 & 451,648.00 & -49.15 & 448,914.00 & -49.46 & 453,820.00 & -48.90 & 447,707.00 & -49.59\\
P5 & 22,689.10 & -56.15 & 19,529.80 & -62.25 & 19,471.70 & -62.37 & 19,000.70 & -63.28 & 18,750.40 & -63.76 & 18,343.40 & -64.55\\
P6 & 508,391.00 & -43.73 & 536,817.00 & -40.58 & 491,931.00 & -45.55 & 481,842.00 & -46.67 & 480,282.00 & -46.84 & 472,967.00 & -47.65\\
P7 & 421,610.00 & -33.73 & 421,067.00 & -33.82 & 417,878.00 & -34.32 & 410,308.00 & -35.51 & 417,676.00 & -34.35 & 415,257.00 & -34.73\\
P8 & 2,126,430.00 & -7.61 & 2,022,190.00 & -12.14 & 1,907,080.00 & -17.14 & 1,782,580.00 & -22.55 & 1,767,330.00 & -23.21 & 1,821,650.00 & -20.85\\ 
\hline
&    & -40.95 &    & -42.13 &    & -43.56 &    & -44.66 &    & -44.65 &    & -45.33\\
\hline
\multicolumn{13}{c}{}\\
\hline
\multicolumn{13}{|c|}{Strategy S9 ($K = 8$)}\\
\hline
\multirow{3}{*}{Inst.} &  
\multicolumn{2}{c|}{600s} &  
\multicolumn{2}{c|}{1200s} & 
\multicolumn{2}{c|}{1800s} & 
\multicolumn{2}{c|}{2400s} & 
\multicolumn{2}{c|}{3000s} & 
\multicolumn{2}{c|}{3600s}\\ 
\cline{2-13}
& \multicolumn{1}{c}{Objective} & \multicolumn{1}{c|}{\multirow{2}{*}{gap (\%)}} 
& \multicolumn{1}{c}{Objective} & \multicolumn{1}{c|}{\multirow{2}{*}{gap (\%)}} 
& \multicolumn{1}{c}{Objective} & \multicolumn{1}{c|}{\multirow{2}{*}{gap (\%)}} 
& \multicolumn{1}{c}{Objective} & \multicolumn{1}{c|}{\multirow{2}{*}{gap (\%)}} 
& \multicolumn{1}{c}{Objective} & \multicolumn{1}{c|}{\multirow{2}{*}{gap (\%)}} 
& \multicolumn{1}{c}{Objective} & \multicolumn{1}{c|}{\multirow{2}{*}{gap (\%)}} \\
& \multicolumn{1}{c}{function value} &
& \multicolumn{1}{c}{function value} &
& \multicolumn{1}{c}{function value} &
& \multicolumn{1}{c}{function value} &
& \multicolumn{1}{c}{function value} &
& \multicolumn{1}{c}{function value} &\\
\hline
\hline
P1 & 641,810.00 & -39.99 & 639,873.00 & -40.17 & 634,966.00 & -40.63 & 634,966.00 & -40.63 & 634,966.00 & -40.63 & 634,966.00 & -40.63\\
P2  & 21,903.40 & -66.15 & 20,775.60 & -67.89 & 20,887.30 & -67.72 & 20,967.60 & -67.60 & 21,158.40 & -67.30 & 20,907.10 & -67.69\\
P3 & 475,705.00 & -36.99 & 475,496.00 & -37.02 & 475,249.00 & -37.05 & 474,957.00 & -37.09 & 475,546.00 & -37.01 & 475,218.00 & -37.05\\
P4 & 468,145.00 & -47.29 & 459,546.00 & -48.26 & 457,263.00 & -48.52 & 461,697.00 & -48.02 & 452,652.00 & -49.04 & 452,628.00 & -49.04\\
P5 & 26,154.90 & -49.45 & 19,585.20 & -62.15 & 20,946.40 & -59.52 & 19,461.40 & -62.39 & 19,163.20 & -62.96 & 18,822.10 & -63.62\\
P6 & 516,010.00 & -42.89 & 478,088.00 & -47.08 & 487,290.00 & -46.07 & 488,026.00 & -45.99 & 475,287.00 & -47.39 & 465,060.00 & -48.53\\
P7 & 424,808.00 & -33.23 & 418,944.00 & -34.15 & 409,001.00 & -35.71 & 410,275.00 & -35.51 & 410,252.00 & -35.52 & 407,847.00 & -35.89\\
P8 & 1,737,070.00 & -24.53 & 1,639,660.00 & -28.76 & 1,666,690.00 & -27.58 & 1,652,060.00 & -28.22 & 1,678,470.00 & -27.07 & 1,543,030.00 & -32.96\\ 
\hline
&    & -42.57 &    & -45.69 &    & -45.35 &    & -45.68 &    & -45.87 &    & -46.93\\
\hline
\end{tabular}}
\end{center}
\caption{Gap between solutions found by CPLEX and strategies S1 and S9 (with $K=8$) and the company's solutions, with varying CPU time limits in $\{600, 1200, 1800, 2400, 3000, 3600\}$.}
\label{tab6}
\end{table}

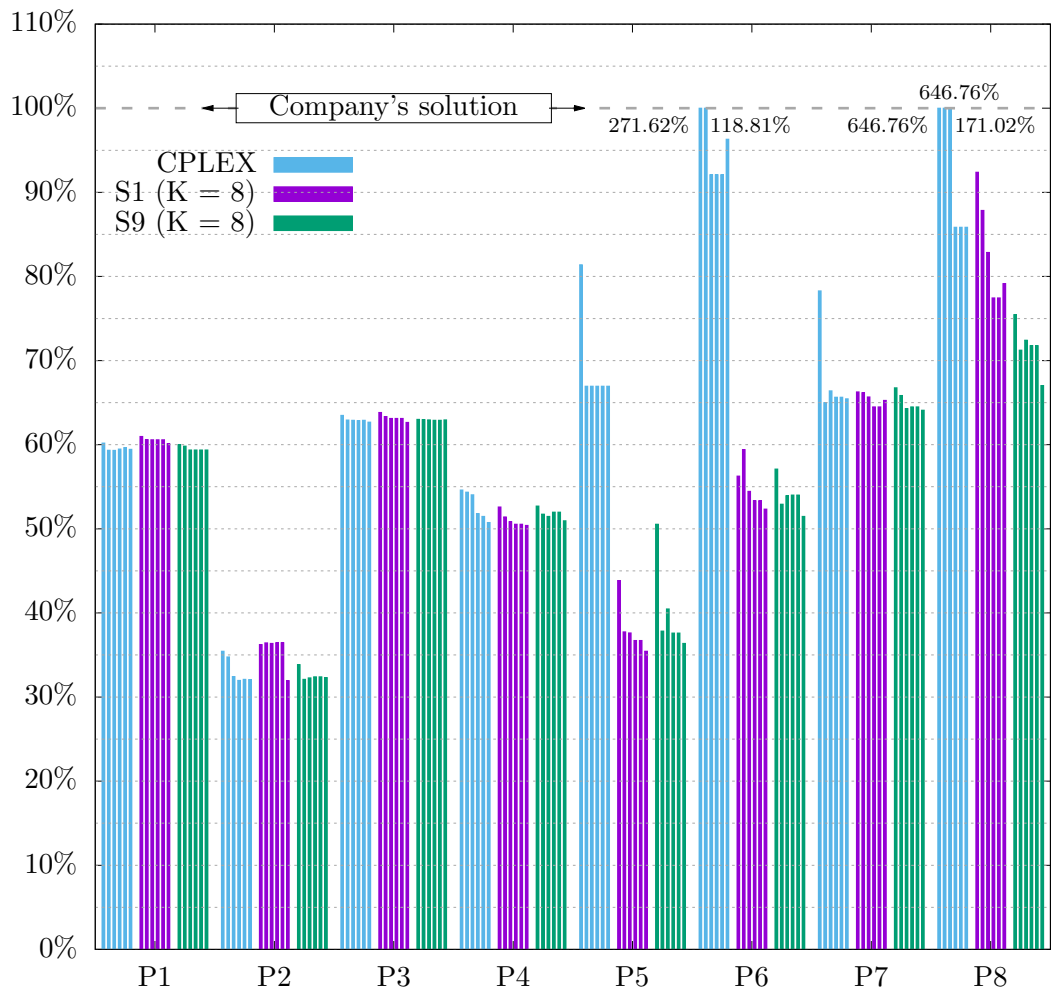
\begin{figure}
\centering
\input{abkm2020fig8.tex}
\caption{Solutions obtained with CPLEX and strategies S1 and S9 (with $S=8$) as a percentage of the company's solutions.}
\label{fig8}
\end{figure}

Figure~\ref{fig8} shows the data from Table~\ref{tab6} in a different way. In the figure, for each instance, the solution reported by the company corresponds to 100\% and the solutions found by the other methods (CPLEX and strategies S1 and S9) appear as percentages of the solution reported by the company. In the figure, CPLEX appears in light blue, strategy~S1 appears in violet and strategy~S9 appears in green. Also, for each method and each instance, there are~6 thin bars that, from left to right, correspond to the time limits $600, 1200, \dots, 3600$, respectively. Bars that exceed 100\% appear truncated in the figure for presentation purposes, and their true values are indicated with numbers near the top. Note, for example, that in instances~P1 and~P3 the~6 thin bars of each method have practically the same height. This means that these instances are apparently simple for the methods, which practically find the same solutions regardless of the time limit. For the other instances, the~6 thin bars associated with different time limits seem to have a decreasing height from left to right, which shows that, with more time, the methods are able to find better solutions.

Overall, as already stated by looking at Table~\ref{tab6}, both CPLEX and strategies~S1 and~S9 improve the solutions reported by the company, except for instances~P6 and~P8 where CPLEX found worse solutions for threshold times below~1800 seconds. Regardless, the performance of CPLEX and strategies~S1 and~S9 are similar in instances~P1, P2, P3 and~P4. In these~4 instances, the~3 methods improve the solution presented by the company and return very similar solutions. It is worth noting that these~4 instances are the smallest instances, with dimensions similar to the random instances of groups~A and~B. In these instances, using CPLEX would be a reasonable alternative. The situation is different in instances P5, P6, P7 and P8, which are the largest within the set considered, both~S1 and~S9 substantially improve the solutions found by the company, in a way that is not being matched by CPLEX. Note even that in instances~P5, P6 and~P8, the solution found by strategies~S1 and~S9 with the shortest time limit is better than the solution found by CPLEX with the longest time limit considered.

\section{Conclusions}

This paper addressed an integrated real-world lot sizing and scheduling problem in a complex operating environment that occurs in a large company in the personal care consumer goods industry. This problem is composed of distinct parallel machines with limited production capacity and sequence-dependent setup times and costs. There is also limited storage capacity for finished goods. In order to solve this problem, a MILP model was presented and several problem-dependent and problem-independent strategies based on the relax-and-fix heuristic were developed. The performance of the heuristics was evaluated by solving randomly generated instances and real-world cases. Exhaustive details about the parameters of the real instances are given (see Table~\ref{tab1}) in such a way that new random instances can be generated that preserve the characteristics of the real instances. This allows the generation of new test sets to evaluate and compare methods that apply to the problem under consideration.

The relax-and-fix strategies introduced showed the best results overall. Strategies S1 and S9 that prioritize chronological time and critical machines, respectively, showed robustly superior performance compared to the other strategies. In addition, the performance of the relax-and-fix strategies was compared with actual results performed by the company and the results obtained by solving the MILP model using CPLEX. Compared to the actual solutions used by the company, the strategies produced results with lower costs and a reduction between 41\% and 47\%, mainly from the reduction of inventory levels and better allocation of production lots on the machines. The relax-and-fix heuristics also outperformed CPLEX applied to the MILP model. Due to the successful application of the proposed methods, the ideas presented in this work are now being extended to more complex production and logistics environments with multiple distribution centers and factories that are also presented in the target company.\\

\noindent
\textbf{Conflict of interest statement:} On behalf of all authors, the corresponding author states that there is no conflict of interest.\\

\noindent
\textbf{Data availability:} The datasets generated during and/or analysed during the current study are available in the GitHub repository, \url{https://github.com/kennedy94/GLSPPL-RF}.

\bibliographystyle{apalike}

\bibliography{abkm2020}

\end{document}

%% file: abkm2020fig3.tex
\begingroup
  \makeatletter
  \providecommand\color[2][]{%
    \GenericError{(gnuplot) \space\space\space\@spaces}{%
      Package color not loaded in conjunction with
      terminal option `colourtext'%
    }{See the gnuplot documentation for explanation.%
    }{Either use 'blacktext' in gnuplot or load the package
      color.sty in LaTeX.}%
    \renewcommand\color[2][]{}%
  }%
  \providecommand\includegraphics[2][]{%
    \GenericError{(gnuplot) \space\space\space\@spaces}{%
      Package graphicx or graphics not loaded%
    }{See the gnuplot documentation for explanation.%
    }{The gnuplot epslatex terminal needs graphicx.sty or graphics.sty.}%
    \renewcommand\includegraphics[2][]{}%
  }%
  \providecommand\rotatebox[2]{#2}%
  \@ifundefined{ifGPcolor}{%
    \newif\ifGPcolor
    \GPcolorfalse
  }{}%
  \@ifundefined{ifGPblacktext}{%
    \newif\ifGPblacktext
    \GPblacktexttrue
  }{}%
  \let\gplgaddtomacro\g@addto@macro
  \gdef\gplbacktext{}%
  \gdef\gplfronttext{}%
  \makeatother
  \ifGPblacktext
    \def\colorrgb#1{}%
    \def\colorgray#1{}%
  \else
    \ifGPcolor
      \def\colorrgb#1{\color[rgb]{#1}}%
      \def\colorgray#1{\color[gray]{#1}}%
      \expandafter\def\csname LTw\endcsname{\color{white}}%
      \expandafter\def\csname LTb\endcsname{\color{black}}%
      \expandafter\def\csname LTa\endcsname{\color{black}}%
      \expandafter\def\csname LT0\endcsname{\color[rgb]{1,0,0}}%
      \expandafter\def\csname LT1\endcsname{\color[rgb]{0,1,0}}%
      \expandafter\def\csname LT2\endcsname{\color[rgb]{0,0,1}}%
      \expandafter\def\csname LT3\endcsname{\color[rgb]{1,0,1}}%
      \expandafter\def\csname LT4\endcsname{\color[rgb]{0,1,1}}%
      \expandafter\def\csname LT5\endcsname{\color[rgb]{1,1,0}}%
      \expandafter\def\csname LT6\endcsname{\color[rgb]{0,0,0}}%
      \expandafter\def\csname LT7\endcsname{\color[rgb]{1,0.3,0}}%
      \expandafter\def\csname LT8\endcsname{\color[rgb]{0.5,0.5,0.5}}%
    \else
      \def\colorrgb#1{\color{black}}%
      \def\colorgray#1{\color[gray]{#1}}%
      \expandafter\def\csname LTw\endcsname{\color{white}}%
      \expandafter\def\csname LTb\endcsname{\color{black}}%
      \expandafter\def\csname LTa\endcsname{\color{black}}%
      \expandafter\def\csname LT0\endcsname{\color{black}}%
      \expandafter\def\csname LT1\endcsname{\color{black}}%
      \expandafter\def\csname LT2\endcsname{\color{black}}%
      \expandafter\def\csname LT3\endcsname{\color{black}}%
      \expandafter\def\csname LT4\endcsname{\color{black}}%
      \expandafter\def\csname LT5\endcsname{\color{black}}%
      \expandafter\def\csname LT6\endcsname{\color{black}}%
      \expandafter\def\csname LT7\endcsname{\color{black}}%
      \expandafter\def\csname LT8\endcsname{\color{black}}%
    \fi
  \fi
    \setlength{\unitlength}{0.0500bp}%
    \ifx\gptboxheight\undefined%
      \newlength{\gptboxheight}%
      \newlength{\gptboxwidth}%
      \newsavebox{\gptboxtext}%
    \fi%
    \setlength{\fboxrule}{0.5pt}%
    \setlength{\fboxsep}{1pt}%
\begin{picture}(6802.00,6802.00)%
    \gplgaddtomacro\gplbacktext{%
      \csname LTb\endcsname
      \put(-132,6345){\makebox(0,0)[r]{\strut{}$0\%$}}%
    }%
    \gplgaddtomacro\gplfronttext{%
      \csname LTb\endcsname
      \put(6966,6231){\rotatebox{-270}{\makebox(0,0){\strut{}Group A}}}%
    }%
    \gplgaddtomacro\gplbacktext{%
      \csname LTb\endcsname
      \put(-132,4681){\makebox(0,0)[r]{\strut{}$0\%$}}%
    }%
    \gplgaddtomacro\gplfronttext{%
      \csname LTb\endcsname
      \put(6966,4871){\rotatebox{-270}{\makebox(0,0){\strut{}Group B}}}%
    }%
    \gplgaddtomacro\gplbacktext{%
      \csname LTb\endcsname
      \put(-132,3977){\makebox(0,0)[r]{\strut{}$0\%$}}%
    }%
    \gplgaddtomacro\gplfronttext{%
      \csname LTb\endcsname
      \put(6966,3510){\rotatebox{-270}{\makebox(0,0){\strut{}Group C}}}%
    }%
    \gplgaddtomacro\gplbacktext{%
      \csname LTb\endcsname
      \put(-132,2625){\makebox(0,0)[r]{\strut{}$0\%$}}%
    }%
    \gplgaddtomacro\gplfronttext{%
      \csname LTb\endcsname
      \put(6966,2150){\rotatebox{-270}{\makebox(0,0){\strut{}Group D}}}%
    }%
    \gplgaddtomacro\gplbacktext{%
      \csname LTb\endcsname
      \put(-132,933){\makebox(0,0)[r]{\strut{}$0\%$}}%
      \put(0,0){\makebox(0,0){\strut{}1}}%
      \put(756,0){\makebox(0,0){\strut{}2}}%
      \put(1511,0){\makebox(0,0){\strut{}3}}%
      \put(2267,0){\makebox(0,0){\strut{}4}}%
      \put(3023,0){\makebox(0,0){\strut{}5}}%
      \put(3778,0){\makebox(0,0){\strut{}6}}%
      \put(4534,0){\makebox(0,0){\strut{}7}}%
      \put(5290,0){\makebox(0,0){\strut{}8}}%
      \put(6045,0){\makebox(0,0){\strut{}9}}%
      \put(6801,0){\makebox(0,0){\strut{}10}}%
      \csname LTb\endcsname
      \put(-943,1645){\rotatebox{90}{\makebox(0,0)[l]{\strut{}Gap to CPLEX solution (i.e., to $K=1$)}}}%
    }%
    \gplgaddtomacro\gplfronttext{%
      \csname LTb\endcsname
      \put(6966,790){\rotatebox{-270}{\makebox(0,0){\strut{}Group E}}}%
    }%
    \gplbacktext
    \put(0,0){\includegraphics{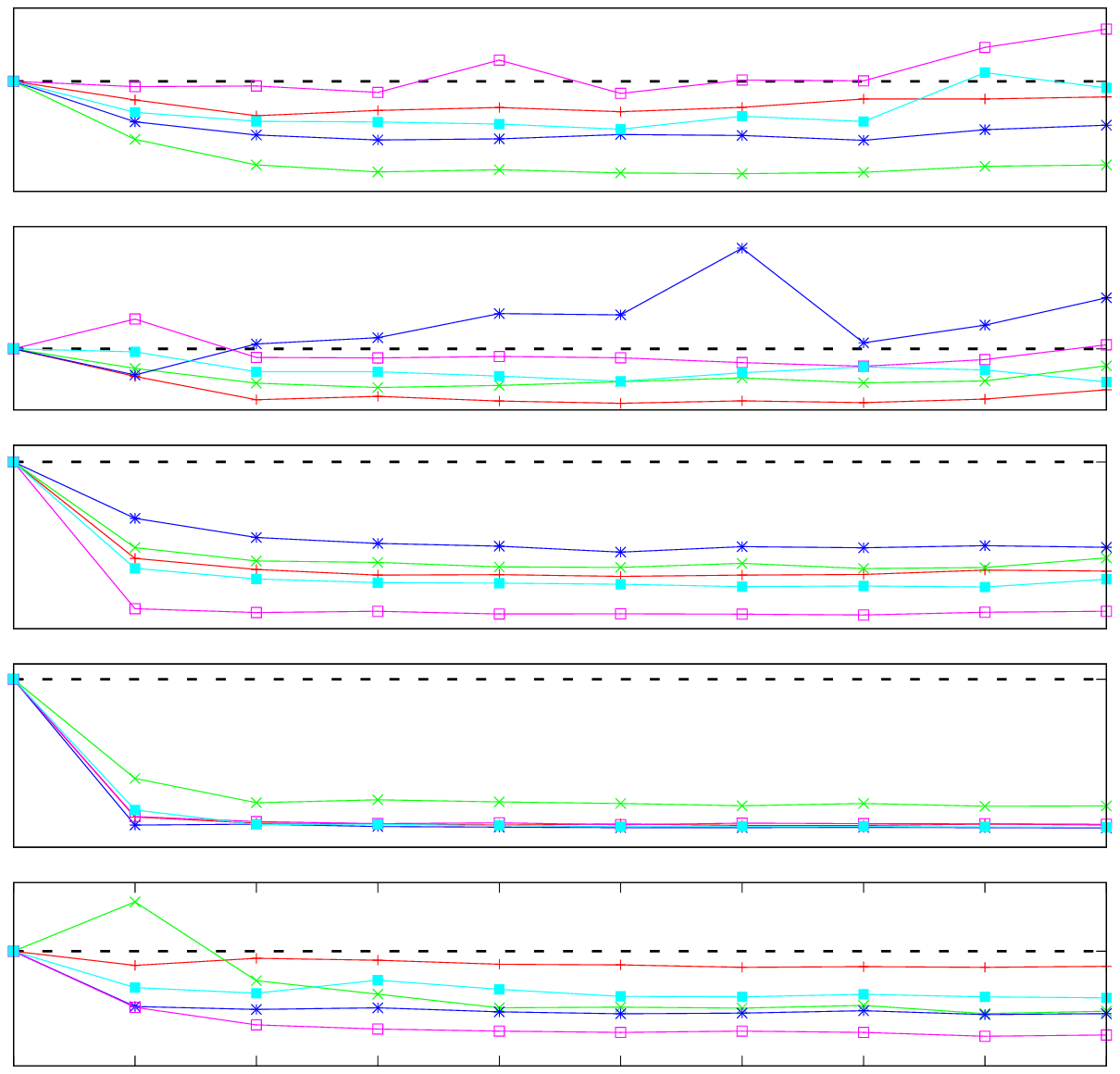}}%
    \gplfronttext
  \end{picture}%
\endgroup

%% file: abkm2020fig4.tex
\begingroup
  \makeatletter
  \providecommand\color[2][]{%
    \GenericError{(gnuplot) \space\space\space\@spaces}{%
      Package color not loaded in conjunction with
      terminal option `colourtext'%
    }{See the gnuplot documentation for explanation.%
    }{Either use 'blacktext' in gnuplot or load the package
      color.sty in LaTeX.}%
    \renewcommand\color[2][]{}%
  }%
  \providecommand\includegraphics[2][]{%
    \GenericError{(gnuplot) \space\space\space\@spaces}{%
      Package graphicx or graphics not loaded%
    }{See the gnuplot documentation for explanation.%
    }{The gnuplot epslatex terminal needs graphicx.sty or graphics.sty.}%
    \renewcommand\includegraphics[2][]{}%
  }%
  \providecommand\rotatebox[2]{#2}%
  \@ifundefined{ifGPcolor}{%
    \newif\ifGPcolor
    \GPcolorfalse
  }{}%
  \@ifundefined{ifGPblacktext}{%
    \newif\ifGPblacktext
    \GPblacktexttrue
  }{}%
  \let\gplgaddtomacro\g@addto@macro
  \gdef\gplbacktext{}%
  \gdef\gplfronttext{}%
  \makeatother
  \ifGPblacktext
    \def\colorrgb#1{}%
    \def\colorgray#1{}%
  \else
    \ifGPcolor
      \def\colorrgb#1{\color[rgb]{#1}}%
      \def\colorgray#1{\color[gray]{#1}}%
      \expandafter\def\csname LTw\endcsname{\color{white}}%
      \expandafter\def\csname LTb\endcsname{\color{black}}%
      \expandafter\def\csname LTa\endcsname{\color{black}}%
      \expandafter\def\csname LT0\endcsname{\color[rgb]{1,0,0}}%
      \expandafter\def\csname LT1\endcsname{\color[rgb]{0,1,0}}%
      \expandafter\def\csname LT2\endcsname{\color[rgb]{0,0,1}}%
      \expandafter\def\csname LT3\endcsname{\color[rgb]{1,0,1}}%
      \expandafter\def\csname LT4\endcsname{\color[rgb]{0,1,1}}%
      \expandafter\def\csname LT5\endcsname{\color[rgb]{1,1,0}}%
      \expandafter\def\csname LT6\endcsname{\color[rgb]{0,0,0}}%
      \expandafter\def\csname LT7\endcsname{\color[rgb]{1,0.3,0}}%
      \expandafter\def\csname LT8\endcsname{\color[rgb]{0.5,0.5,0.5}}%
    \else
      \def\colorrgb#1{\color{black}}%
      \def\colorgray#1{\color[gray]{#1}}%
      \expandafter\def\csname LTw\endcsname{\color{white}}%
      \expandafter\def\csname LTb\endcsname{\color{black}}%
      \expandafter\def\csname LTa\endcsname{\color{black}}%
      \expandafter\def\csname LT0\endcsname{\color{black}}%
      \expandafter\def\csname LT1\endcsname{\color{black}}%
      \expandafter\def\csname LT2\endcsname{\color{black}}%
      \expandafter\def\csname LT3\endcsname{\color{black}}%
      \expandafter\def\csname LT4\endcsname{\color{black}}%
      \expandafter\def\csname LT5\endcsname{\color{black}}%
      \expandafter\def\csname LT6\endcsname{\color{black}}%
      \expandafter\def\csname LT7\endcsname{\color{black}}%
      \expandafter\def\csname LT8\endcsname{\color{black}}%
    \fi
  \fi
    \setlength{\unitlength}{0.0500bp}%
    \ifx\gptboxheight\undefined%
      \newlength{\gptboxheight}%
      \newlength{\gptboxwidth}%
      \newsavebox{\gptboxtext}%
    \fi%
    \setlength{\fboxrule}{0.5pt}%
    \setlength{\fboxsep}{1pt}%
\begin{picture}(6802.00,6802.00)%
    \gplgaddtomacro\gplbacktext{%
      \csname LTb\endcsname
      \put(-132,6345){\makebox(0,0)[r]{\strut{}$0\%$}}%
    }%
    \gplgaddtomacro\gplfronttext{%
      \csname LTb\endcsname
      \put(6966,6231){\rotatebox{-270}{\makebox(0,0){\strut{}Group A}}}%
    }%
    \gplgaddtomacro\gplbacktext{%
      \csname LTb\endcsname
      \put(-132,5251){\makebox(0,0)[r]{\strut{}$0\%$}}%
    }%
    \gplgaddtomacro\gplfronttext{%
      \csname LTb\endcsname
      \put(6966,4871){\rotatebox{-270}{\makebox(0,0){\strut{}Group B}}}%
    }%
    \gplgaddtomacro\gplbacktext{%
      \csname LTb\endcsname
      \put(-132,3977){\makebox(0,0)[r]{\strut{}$0\%$}}%
    }%
    \gplgaddtomacro\gplfronttext{%
      \csname LTb\endcsname
      \put(6966,3510){\rotatebox{-270}{\makebox(0,0){\strut{}Group C}}}%
    }%
    \gplgaddtomacro\gplbacktext{%
      \csname LTb\endcsname
      \put(-132,2625){\makebox(0,0)[r]{\strut{}$0\%$}}%
    }%
    \gplgaddtomacro\gplfronttext{%
      \csname LTb\endcsname
      \put(6966,2150){\rotatebox{-270}{\makebox(0,0){\strut{}Group D}}}%
    }%
    \gplgaddtomacro\gplbacktext{%
      \csname LTb\endcsname
      \put(-132,933){\makebox(0,0)[r]{\strut{}$0\%$}}%
      \put(0,0){\makebox(0,0){\strut{}1}}%
      \put(756,0){\makebox(0,0){\strut{}2}}%
      \put(1511,0){\makebox(0,0){\strut{}3}}%
      \put(2267,0){\makebox(0,0){\strut{}4}}%
      \put(3023,0){\makebox(0,0){\strut{}5}}%
      \put(3778,0){\makebox(0,0){\strut{}6}}%
      \put(4534,0){\makebox(0,0){\strut{}7}}%
      \put(5290,0){\makebox(0,0){\strut{}8}}%
      \put(6045,0){\makebox(0,0){\strut{}9}}%
      \put(6801,0){\makebox(0,0){\strut{}10}}%
      \csname LTb\endcsname
      \put(-943,1645){\rotatebox{90}{\makebox(0,0)[l]{\strut{}Gap to CPLEX solution (i.e., to $K=1$)}}}%
    }%
    \gplgaddtomacro\gplfronttext{%
      \csname LTb\endcsname
      \put(6966,790){\rotatebox{-270}{\makebox(0,0){\strut{}Group E}}}%
    }%
    \gplbacktext
    \put(0,0){\includegraphics{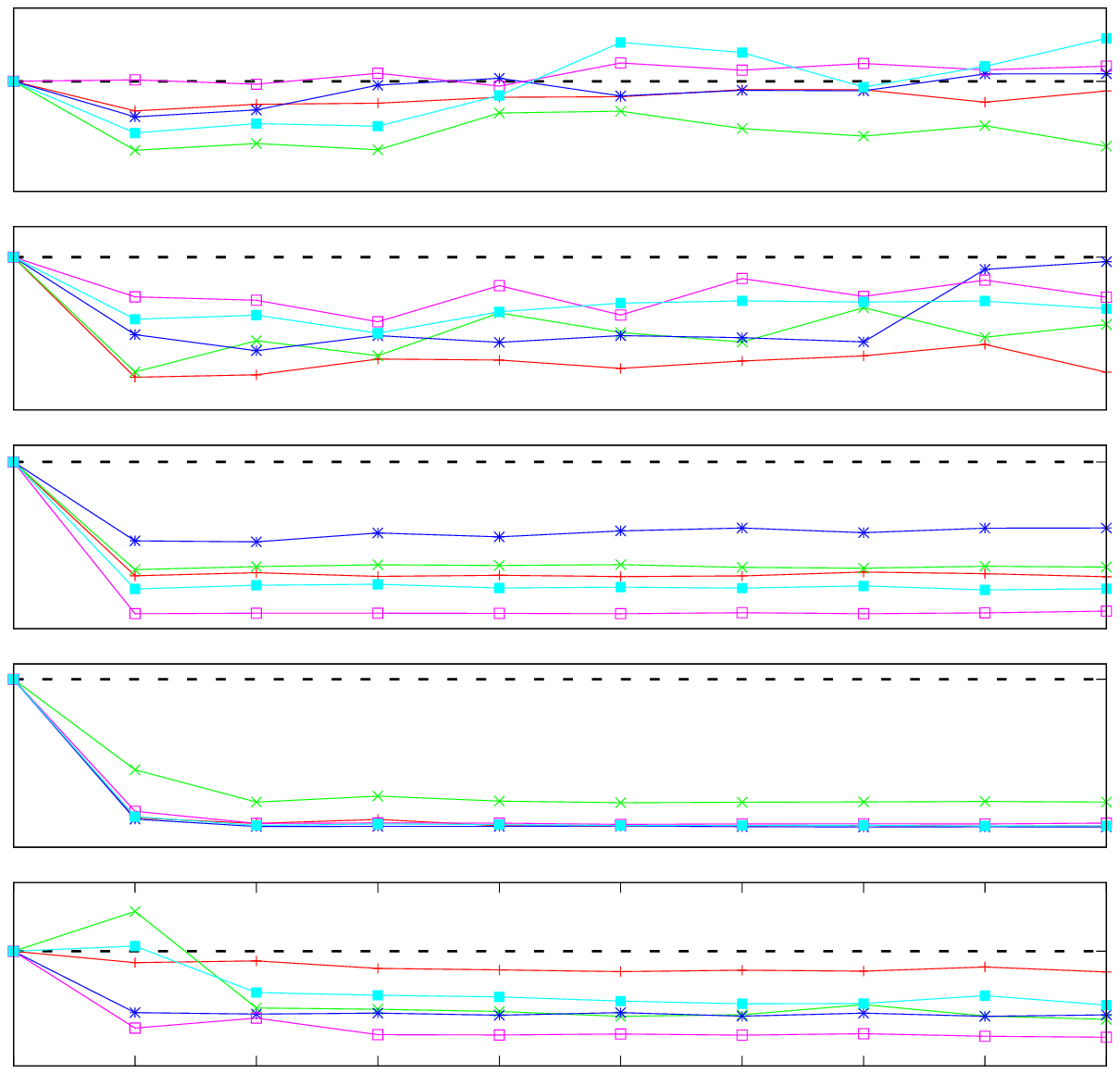}}%
    \gplfronttext
  \end{picture}%
\endgroup

%% file: abkm2020fig5a.tex
\begingroup
  \makeatletter
  \providecommand\color[2][]{%
    \GenericError{(gnuplot) \space\space\space\@spaces}{%
      Package color not loaded in conjunction with
      terminal option `colourtext'%
    }{See the gnuplot documentation for explanation.%
    }{Either use 'blacktext' in gnuplot or load the package
      color.sty in LaTeX.}%
    \renewcommand\color[2][]{}%
  }%
  \providecommand\includegraphics[2][]{%
    \GenericError{(gnuplot) \space\space\space\@spaces}{%
      Package graphicx or graphics not loaded%
    }{See the gnuplot documentation for explanation.%
    }{The gnuplot epslatex terminal needs graphicx.sty or graphics.sty.}%
    \renewcommand\includegraphics[2][]{}%
  }%
  \providecommand\rotatebox[2]{#2}%
  \@ifundefined{ifGPcolor}{%
    \newif\ifGPcolor
    \GPcolorfalse
  }{}%
  \@ifundefined{ifGPblacktext}{%
    \newif\ifGPblacktext
    \GPblacktexttrue
  }{}%
  \let\gplgaddtomacro\g@addto@macro
  \gdef\gplbacktext{}%
  \gdef\gplfronttext{}%
  \makeatother
  \ifGPblacktext
    \def\colorrgb#1{}%
    \def\colorgray#1{}%
  \else
    \ifGPcolor
      \def\colorrgb#1{\color[rgb]{#1}}%
      \def\colorgray#1{\color[gray]{#1}}%
      \expandafter\def\csname LTw\endcsname{\color{white}}%
      \expandafter\def\csname LTb\endcsname{\color{black}}%
      \expandafter\def\csname LTa\endcsname{\color{black}}%
      \expandafter\def\csname LT0\endcsname{\color[rgb]{1,0,0}}%
      \expandafter\def\csname LT1\endcsname{\color[rgb]{0,1,0}}%
      \expandafter\def\csname LT2\endcsname{\color[rgb]{0,0,1}}%
      \expandafter\def\csname LT3\endcsname{\color[rgb]{1,0,1}}%
      \expandafter\def\csname LT4\endcsname{\color[rgb]{0,1,1}}%
      \expandafter\def\csname LT5\endcsname{\color[rgb]{1,1,0}}%
      \expandafter\def\csname LT6\endcsname{\color[rgb]{0,0,0}}%
      \expandafter\def\csname LT7\endcsname{\color[rgb]{1,0.3,0}}%
      \expandafter\def\csname LT8\endcsname{\color[rgb]{0.5,0.5,0.5}}%
    \else
      \def\colorrgb#1{\color{black}}%
      \def\colorgray#1{\color[gray]{#1}}%
      \expandafter\def\csname LTw\endcsname{\color{white}}%
      \expandafter\def\csname LTb\endcsname{\color{black}}%
      \expandafter\def\csname LTa\endcsname{\color{black}}%
      \expandafter\def\csname LT0\endcsname{\color{black}}%
      \expandafter\def\csname LT1\endcsname{\color{black}}%
      \expandafter\def\csname LT2\endcsname{\color{black}}%
      \expandafter\def\csname LT3\endcsname{\color{black}}%
      \expandafter\def\csname LT4\endcsname{\color{black}}%
      \expandafter\def\csname LT5\endcsname{\color{black}}%
      \expandafter\def\csname LT6\endcsname{\color{black}}%
      \expandafter\def\csname LT7\endcsname{\color{black}}%
      \expandafter\def\csname LT8\endcsname{\color{black}}%
    \fi
  \fi
    \setlength{\unitlength}{0.0500bp}%
    \ifx\gptboxheight\undefined%
      \newlength{\gptboxheight}%
      \newlength{\gptboxwidth}%
      \newsavebox{\gptboxtext}%
    \fi%
    \setlength{\fboxrule}{0.5pt}%
    \setlength{\fboxsep}{1pt}%
\begin{picture}(6802.00,3400.00)%
    \gplgaddtomacro\gplbacktext{%
      \csname LTb\endcsname
      \put(990,440){\makebox(0,0)[r]{\strut{}$-100\%$}}%
      \put(990,831){\makebox(0,0)[r]{\strut{}$-80\%$}}%
      \put(990,1223){\makebox(0,0)[r]{\strut{}$-60\%$}}%
      \put(990,1614){\makebox(0,0)[r]{\strut{}$-40\%$}}%
      \put(990,2005){\makebox(0,0)[r]{\strut{}$-20\%$}}%
      \put(990,2396){\makebox(0,0)[r]{\strut{}$0\%$}}%
      \put(990,2788){\makebox(0,0)[r]{\strut{}$20\%$}}%
      \put(990,3179){\makebox(0,0)[r]{\strut{}$40\%$}}%
      \put(1732,220){\makebox(0,0){\strut{}600s}}%
      \put(2544,220){\makebox(0,0){\strut{}1200s}}%
      \put(3357,220){\makebox(0,0){\strut{}1800s}}%
      \put(4170,220){\makebox(0,0){\strut{}2400s}}%
      \put(4983,220){\makebox(0,0){\strut{}3000s}}%
      \put(5795,220){\makebox(0,0){\strut{}3600s}}%
    }%
    \gplgaddtomacro\gplfronttext{%
    }%
    \gplbacktext
    \put(0,0){\includegraphics{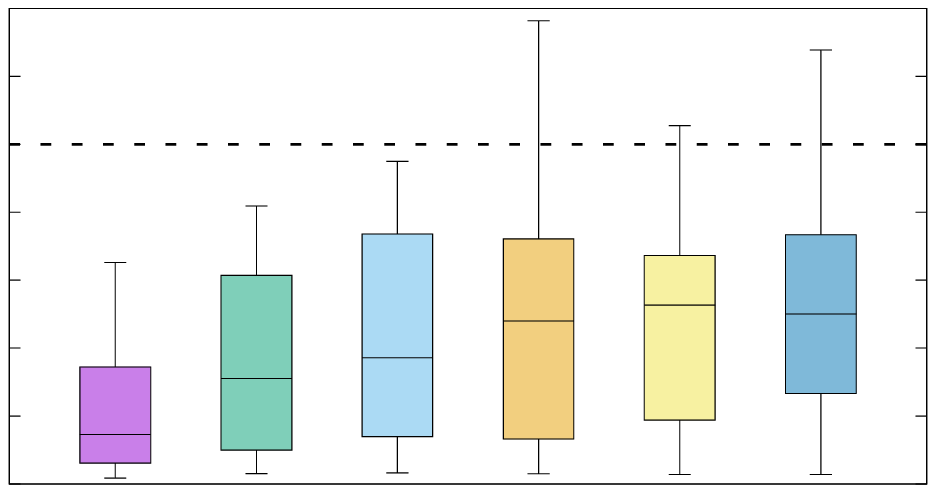}}%
    \gplfronttext
  \end{picture}%
\endgroup

%% file: abkm2020fig5b.tex
\begingroup
  \makeatletter
  \providecommand\color[2][]{%
    \GenericError{(gnuplot) \space\space\space\@spaces}{%
      Package color not loaded in conjunction with
      terminal option `colourtext'%
    }{See the gnuplot documentation for explanation.%
    }{Either use 'blacktext' in gnuplot or load the package
      color.sty in LaTeX.}%
    \renewcommand\color[2][]{}%
  }%
  \providecommand\includegraphics[2][]{%
    \GenericError{(gnuplot) \space\space\space\@spaces}{%
      Package graphicx or graphics not loaded%
    }{See the gnuplot documentation for explanation.%
    }{The gnuplot epslatex terminal needs graphicx.sty or graphics.sty.}%
    \renewcommand\includegraphics[2][]{}%
  }%
  \providecommand\rotatebox[2]{#2}%
  \@ifundefined{ifGPcolor}{%
    \newif\ifGPcolor
    \GPcolorfalse
  }{}%
  \@ifundefined{ifGPblacktext}{%
    \newif\ifGPblacktext
    \GPblacktexttrue
  }{}%
  \let\gplgaddtomacro\g@addto@macro
  \gdef\gplbacktext{}%
  \gdef\gplfronttext{}%
  \makeatother
  \ifGPblacktext
    \def\colorrgb#1{}%
    \def\colorgray#1{}%
  \else
    \ifGPcolor
      \def\colorrgb#1{\color[rgb]{#1}}%
      \def\colorgray#1{\color[gray]{#1}}%
      \expandafter\def\csname LTw\endcsname{\color{white}}%
      \expandafter\def\csname LTb\endcsname{\color{black}}%
      \expandafter\def\csname LTa\endcsname{\color{black}}%
      \expandafter\def\csname LT0\endcsname{\color[rgb]{1,0,0}}%
      \expandafter\def\csname LT1\endcsname{\color[rgb]{0,1,0}}%
      \expandafter\def\csname LT2\endcsname{\color[rgb]{0,0,1}}%
      \expandafter\def\csname LT3\endcsname{\color[rgb]{1,0,1}}%
      \expandafter\def\csname LT4\endcsname{\color[rgb]{0,1,1}}%
      \expandafter\def\csname LT5\endcsname{\color[rgb]{1,1,0}}%
      \expandafter\def\csname LT6\endcsname{\color[rgb]{0,0,0}}%
      \expandafter\def\csname LT7\endcsname{\color[rgb]{1,0.3,0}}%
      \expandafter\def\csname LT8\endcsname{\color[rgb]{0.5,0.5,0.5}}%
    \else
      \def\colorrgb#1{\color{black}}%
      \def\colorgray#1{\color[gray]{#1}}%
      \expandafter\def\csname LTw\endcsname{\color{white}}%
      \expandafter\def\csname LTb\endcsname{\color{black}}%
      \expandafter\def\csname LTa\endcsname{\color{black}}%
      \expandafter\def\csname LT0\endcsname{\color{black}}%
      \expandafter\def\csname LT1\endcsname{\color{black}}%
      \expandafter\def\csname LT2\endcsname{\color{black}}%
      \expandafter\def\csname LT3\endcsname{\color{black}}%
      \expandafter\def\csname LT4\endcsname{\color{black}}%
      \expandafter\def\csname LT5\endcsname{\color{black}}%
      \expandafter\def\csname LT6\endcsname{\color{black}}%
      \expandafter\def\csname LT7\endcsname{\color{black}}%
      \expandafter\def\csname LT8\endcsname{\color{black}}%
    \fi
  \fi
    \setlength{\unitlength}{0.0500bp}%
    \ifx\gptboxheight\undefined%
      \newlength{\gptboxheight}%
      \newlength{\gptboxwidth}%
      \newsavebox{\gptboxtext}%
    \fi%
    \setlength{\fboxrule}{0.5pt}%
    \setlength{\fboxsep}{1pt}%
\begin{picture}(6802.00,3400.00)%
    \gplgaddtomacro\gplbacktext{%
      \csname LTb\endcsname
      \put(726,744){\makebox(0,0)[r]{\strut{}$0\%$}}%
      \put(726,1353){\makebox(0,0)[r]{\strut{}$20\%$}}%
      \put(726,1962){\makebox(0,0)[r]{\strut{}$40\%$}}%
      \put(726,2570){\makebox(0,0)[r]{\strut{}$60\%$}}%
      \put(726,3179){\makebox(0,0)[r]{\strut{}$80\%$}}%
      \put(1614,220){\makebox(0,0){\strut{}600s}}%
      \put(2623,220){\makebox(0,0){\strut{}1200s}}%
      \put(3632,220){\makebox(0,0){\strut{}1800s}}%
      \put(4640,220){\makebox(0,0){\strut{}2400s}}%
      \put(5649,220){\makebox(0,0){\strut{}3000s}}%
    }%
    \gplgaddtomacro\gplfronttext{%
    }%
    \gplbacktext
    \put(0,0){\includegraphics{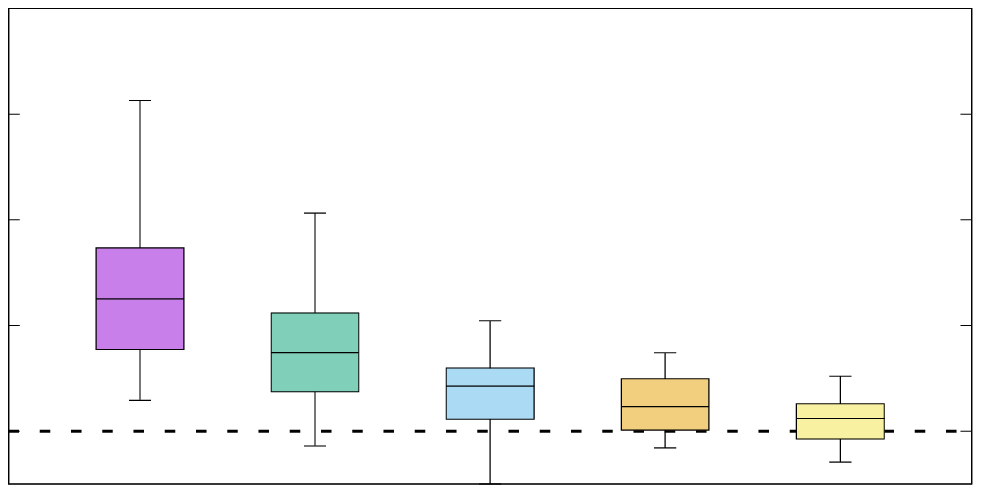}}%
    \gplfronttext
  \end{picture}%
\endgroup

%% file: abkm2020fig6a.tex
\begingroup
  \makeatletter
  \providecommand\color[2][]{%
    \GenericError{(gnuplot) \space\space\space\@spaces}{%
      Package color not loaded in conjunction with
      terminal option `colourtext'%
    }{See the gnuplot documentation for explanation.%
    }{Either use 'blacktext' in gnuplot or load the package
      color.sty in LaTeX.}%
    \renewcommand\color[2][]{}%
  }%
  \providecommand\includegraphics[2][]{%
    \GenericError{(gnuplot) \space\space\space\@spaces}{%
      Package graphicx or graphics not loaded%
    }{See the gnuplot documentation for explanation.%
    }{The gnuplot epslatex terminal needs graphicx.sty or graphics.sty.}%
    \renewcommand\includegraphics[2][]{}%
  }%
  \providecommand\rotatebox[2]{#2}%
  \@ifundefined{ifGPcolor}{%
    \newif\ifGPcolor
    \GPcolorfalse
  }{}%
  \@ifundefined{ifGPblacktext}{%
    \newif\ifGPblacktext
    \GPblacktexttrue
  }{}%
  \let\gplgaddtomacro\g@addto@macro
  \gdef\gplbacktext{}%
  \gdef\gplfronttext{}%
  \makeatother
  \ifGPblacktext
    \def\colorrgb#1{}%
    \def\colorgray#1{}%
  \else
    \ifGPcolor
      \def\colorrgb#1{\color[rgb]{#1}}%
      \def\colorgray#1{\color[gray]{#1}}%
      \expandafter\def\csname LTw\endcsname{\color{white}}%
      \expandafter\def\csname LTb\endcsname{\color{black}}%
      \expandafter\def\csname LTa\endcsname{\color{black}}%
      \expandafter\def\csname LT0\endcsname{\color[rgb]{1,0,0}}%
      \expandafter\def\csname LT1\endcsname{\color[rgb]{0,1,0}}%
      \expandafter\def\csname LT2\endcsname{\color[rgb]{0,0,1}}%
      \expandafter\def\csname LT3\endcsname{\color[rgb]{1,0,1}}%
      \expandafter\def\csname LT4\endcsname{\color[rgb]{0,1,1}}%
      \expandafter\def\csname LT5\endcsname{\color[rgb]{1,1,0}}%
      \expandafter\def\csname LT6\endcsname{\color[rgb]{0,0,0}}%
      \expandafter\def\csname LT7\endcsname{\color[rgb]{1,0.3,0}}%
      \expandafter\def\csname LT8\endcsname{\color[rgb]{0.5,0.5,0.5}}%
    \else
      \def\colorrgb#1{\color{black}}%
      \def\colorgray#1{\color[gray]{#1}}%
      \expandafter\def\csname LTw\endcsname{\color{white}}%
      \expandafter\def\csname LTb\endcsname{\color{black}}%
      \expandafter\def\csname LTa\endcsname{\color{black}}%
      \expandafter\def\csname LT0\endcsname{\color{black}}%
      \expandafter\def\csname LT1\endcsname{\color{black}}%
      \expandafter\def\csname LT2\endcsname{\color{black}}%
      \expandafter\def\csname LT3\endcsname{\color{black}}%
      \expandafter\def\csname LT4\endcsname{\color{black}}%
      \expandafter\def\csname LT5\endcsname{\color{black}}%
      \expandafter\def\csname LT6\endcsname{\color{black}}%
      \expandafter\def\csname LT7\endcsname{\color{black}}%
      \expandafter\def\csname LT8\endcsname{\color{black}}%
    \fi
  \fi
    \setlength{\unitlength}{0.0500bp}%
    \ifx\gptboxheight\undefined%
      \newlength{\gptboxheight}%
      \newlength{\gptboxwidth}%
      \newsavebox{\gptboxtext}%
    \fi%
    \setlength{\fboxrule}{0.5pt}%
    \setlength{\fboxsep}{1pt}%
\begin{picture}(6802.00,3400.00)%
    \gplgaddtomacro\gplbacktext{%
      \csname LTb\endcsname
      \put(990,440){\makebox(0,0)[r]{\strut{}$-100\%$}}%
      \put(990,938){\makebox(0,0)[r]{\strut{}$-80\%$}}%
      \put(990,1436){\makebox(0,0)[r]{\strut{}$-60\%$}}%
      \put(990,1934){\makebox(0,0)[r]{\strut{}$-40\%$}}%
      \put(990,2432){\makebox(0,0)[r]{\strut{}$-20\%$}}%
      \put(990,2930){\makebox(0,0)[r]{\strut{}$0\%$}}%
      \put(1732,220){\makebox(0,0){\strut{}600s}}%
      \put(2544,220){\makebox(0,0){\strut{}1200s}}%
      \put(3357,220){\makebox(0,0){\strut{}1800s}}%
      \put(4170,220){\makebox(0,0){\strut{}2400s}}%
      \put(4983,220){\makebox(0,0){\strut{}3000s}}%
      \put(5795,220){\makebox(0,0){\strut{}3600s}}%
    }%
    \gplgaddtomacro\gplfronttext{%
    }%
    \gplbacktext
    \put(0,0){\includegraphics{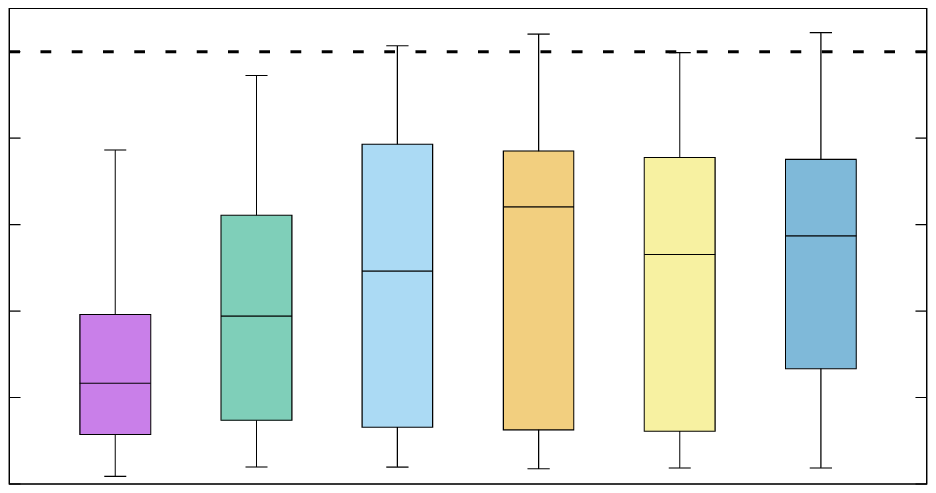}}%
    \gplfronttext
  \end{picture}%
\endgroup

%% file: abkm2020fig6b.tex
\begingroup
  \makeatletter
  \providecommand\color[2][]{%
    \GenericError{(gnuplot) \space\space\space\@spaces}{%
      Package color not loaded in conjunction with
      terminal option `colourtext'%
    }{See the gnuplot documentation for explanation.%
    }{Either use 'blacktext' in gnuplot or load the package
      color.sty in LaTeX.}%
    \renewcommand\color[2][]{}%
  }%
  \providecommand\includegraphics[2][]{%
    \GenericError{(gnuplot) \space\space\space\@spaces}{%
      Package graphicx or graphics not loaded%
    }{See the gnuplot documentation for explanation.%
    }{The gnuplot epslatex terminal needs graphicx.sty or graphics.sty.}%
    \renewcommand\includegraphics[2][]{}%
  }%
  \providecommand\rotatebox[2]{#2}%
  \@ifundefined{ifGPcolor}{%
    \newif\ifGPcolor
    \GPcolorfalse
  }{}%
  \@ifundefined{ifGPblacktext}{%
    \newif\ifGPblacktext
    \GPblacktexttrue
  }{}%
  \let\gplgaddtomacro\g@addto@macro
  \gdef\gplbacktext{}%
  \gdef\gplfronttext{}%
  \makeatother
  \ifGPblacktext
    \def\colorrgb#1{}%
    \def\colorgray#1{}%
  \else
    \ifGPcolor
      \def\colorrgb#1{\color[rgb]{#1}}%
      \def\colorgray#1{\color[gray]{#1}}%
      \expandafter\def\csname LTw\endcsname{\color{white}}%
      \expandafter\def\csname LTb\endcsname{\color{black}}%
      \expandafter\def\csname LTa\endcsname{\color{black}}%
      \expandafter\def\csname LT0\endcsname{\color[rgb]{1,0,0}}%
      \expandafter\def\csname LT1\endcsname{\color[rgb]{0,1,0}}%
      \expandafter\def\csname LT2\endcsname{\color[rgb]{0,0,1}}%
      \expandafter\def\csname LT3\endcsname{\color[rgb]{1,0,1}}%
      \expandafter\def\csname LT4\endcsname{\color[rgb]{0,1,1}}%
      \expandafter\def\csname LT5\endcsname{\color[rgb]{1,1,0}}%
      \expandafter\def\csname LT6\endcsname{\color[rgb]{0,0,0}}%
      \expandafter\def\csname LT7\endcsname{\color[rgb]{1,0.3,0}}%
      \expandafter\def\csname LT8\endcsname{\color[rgb]{0.5,0.5,0.5}}%
    \else
      \def\colorrgb#1{\color{black}}%
      \def\colorgray#1{\color[gray]{#1}}%
      \expandafter\def\csname LTw\endcsname{\color{white}}%
      \expandafter\def\csname LTb\endcsname{\color{black}}%
      \expandafter\def\csname LTa\endcsname{\color{black}}%
      \expandafter\def\csname LT0\endcsname{\color{black}}%
      \expandafter\def\csname LT1\endcsname{\color{black}}%
      \expandafter\def\csname LT2\endcsname{\color{black}}%
      \expandafter\def\csname LT3\endcsname{\color{black}}%
      \expandafter\def\csname LT4\endcsname{\color{black}}%
      \expandafter\def\csname LT5\endcsname{\color{black}}%
      \expandafter\def\csname LT6\endcsname{\color{black}}%
      \expandafter\def\csname LT7\endcsname{\color{black}}%
      \expandafter\def\csname LT8\endcsname{\color{black}}%
    \fi
  \fi
    \setlength{\unitlength}{0.0500bp}%
    \ifx\gptboxheight\undefined%
      \newlength{\gptboxheight}%
      \newlength{\gptboxwidth}%
      \newsavebox{\gptboxtext}%
    \fi%
    \setlength{\fboxrule}{0.5pt}%
    \setlength{\fboxsep}{1pt}%
\begin{picture}(6802.00,3400.00)%
    \gplgaddtomacro\gplbacktext{%
      \csname LTb\endcsname
      \put(858,440){\makebox(0,0)[r]{\strut{}$-20\%$}}%
      \put(858,1223){\makebox(0,0)[r]{\strut{}$0\%$}}%
      \put(858,2005){\makebox(0,0)[r]{\strut{}$20\%$}}%
      \put(858,2788){\makebox(0,0)[r]{\strut{}$40\%$}}%
      \put(1728,220){\makebox(0,0){\strut{}600s}}%
      \put(2713,220){\makebox(0,0){\strut{}1200s}}%
      \put(3698,220){\makebox(0,0){\strut{}1800s}}%
      \put(4682,220){\makebox(0,0){\strut{}2400s}}%
      \put(5667,220){\makebox(0,0){\strut{}3000s}}%
    }%
    \gplgaddtomacro\gplfronttext{%
    }%
    \gplbacktext
    \put(0,0){\includegraphics{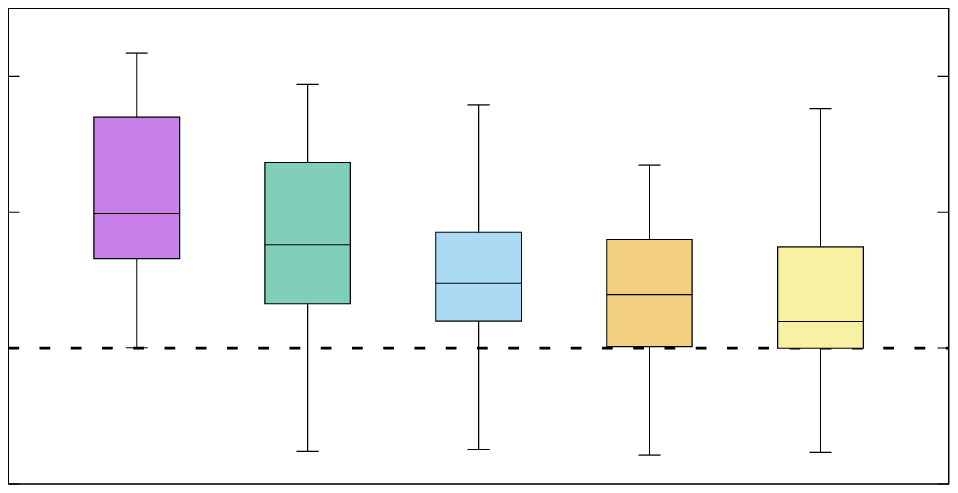}}%
    \gplfronttext
  \end{picture}%
\endgroup

%% file: abkm2020fig7a.tex
\begingroup
  \makeatletter
  \providecommand\color[2][]{%
    \GenericError{(gnuplot) \space\space\space\@spaces}{%
      Package color not loaded in conjunction with
      terminal option `colourtext'%
    }{See the gnuplot documentation for explanation.%
    }{Either use 'blacktext' in gnuplot or load the package
      color.sty in LaTeX.}%
    \renewcommand\color[2][]{}%
  }%
  \providecommand\includegraphics[2][]{%
    \GenericError{(gnuplot) \space\space\space\@spaces}{%
      Package graphicx or graphics not loaded%
    }{See the gnuplot documentation for explanation.%
    }{The gnuplot epslatex terminal needs graphicx.sty or graphics.sty.}%
    \renewcommand\includegraphics[2][]{}%
  }%
  \providecommand\rotatebox[2]{#2}%
  \@ifundefined{ifGPcolor}{%
    \newif\ifGPcolor
    \GPcolorfalse
  }{}%
  \@ifundefined{ifGPblacktext}{%
    \newif\ifGPblacktext
    \GPblacktexttrue
  }{}%
  \let\gplgaddtomacro\g@addto@macro
  \gdef\gplbacktext{}%
  \gdef\gplfronttext{}%
  \makeatother
  \ifGPblacktext
    \def\colorrgb#1{}%
    \def\colorgray#1{}%
  \else
    \ifGPcolor
      \def\colorrgb#1{\color[rgb]{#1}}%
      \def\colorgray#1{\color[gray]{#1}}%
      \expandafter\def\csname LTw\endcsname{\color{white}}%
      \expandafter\def\csname LTb\endcsname{\color{black}}%
      \expandafter\def\csname LTa\endcsname{\color{black}}%
      \expandafter\def\csname LT0\endcsname{\color[rgb]{1,0,0}}%
      \expandafter\def\csname LT1\endcsname{\color[rgb]{0,1,0}}%
      \expandafter\def\csname LT2\endcsname{\color[rgb]{0,0,1}}%
      \expandafter\def\csname LT3\endcsname{\color[rgb]{1,0,1}}%
      \expandafter\def\csname LT4\endcsname{\color[rgb]{0,1,1}}%
      \expandafter\def\csname LT5\endcsname{\color[rgb]{1,1,0}}%
      \expandafter\def\csname LT6\endcsname{\color[rgb]{0,0,0}}%
      \expandafter\def\csname LT7\endcsname{\color[rgb]{1,0.3,0}}%
      \expandafter\def\csname LT8\endcsname{\color[rgb]{0.5,0.5,0.5}}%
    \else
      \def\colorrgb#1{\color{black}}%
      \def\colorgray#1{\color[gray]{#1}}%
      \expandafter\def\csname LTw\endcsname{\color{white}}%
      \expandafter\def\csname LTb\endcsname{\color{black}}%
      \expandafter\def\csname LTa\endcsname{\color{black}}%
      \expandafter\def\csname LT0\endcsname{\color{black}}%
      \expandafter\def\csname LT1\endcsname{\color{black}}%
      \expandafter\def\csname LT2\endcsname{\color{black}}%
      \expandafter\def\csname LT3\endcsname{\color{black}}%
      \expandafter\def\csname LT4\endcsname{\color{black}}%
      \expandafter\def\csname LT5\endcsname{\color{black}}%
      \expandafter\def\csname LT6\endcsname{\color{black}}%
      \expandafter\def\csname LT7\endcsname{\color{black}}%
      \expandafter\def\csname LT8\endcsname{\color{black}}%
    \fi
  \fi
    \setlength{\unitlength}{0.0500bp}%
    \ifx\gptboxheight\undefined%
      \newlength{\gptboxheight}%
      \newlength{\gptboxwidth}%
      \newsavebox{\gptboxtext}%
    \fi%
    \setlength{\fboxrule}{0.5pt}%
    \setlength{\fboxsep}{1pt}%
\begin{picture}(3400.00,3400.00)%
    \gplgaddtomacro\gplbacktext{%
      \csname LTb\endcsname
      \put(-132,220){\makebox(0,0)[r]{\strut{}$-50\%$}}%
      \put(-132,750){\makebox(0,0)[r]{\strut{}$-40\%$}}%
      \put(-132,1280){\makebox(0,0)[r]{\strut{}$-30\%$}}%
      \put(-132,1810){\makebox(0,0)[r]{\strut{}$-20\%$}}%
      \put(-132,2339){\makebox(0,0)[r]{\strut{}$-10\%$}}%
      \put(-132,2869){\makebox(0,0)[r]{\strut{}$0\%$}}%
      \put(-132,3399){\makebox(0,0)[r]{\strut{}$10\%$}}%
      \put(0,0){\makebox(0,0){\strut{}2}}%
      \put(425,0){\makebox(0,0){\strut{}3}}%
      \put(850,0){\makebox(0,0){\strut{}4}}%
      \put(1275,0){\makebox(0,0){\strut{}5}}%
      \put(1700,0){\makebox(0,0){\strut{}6}}%
      \put(2124,0){\makebox(0,0){\strut{}7}}%
      \put(2549,0){\makebox(0,0){\strut{}8}}%
      \put(2974,0){\makebox(0,0){\strut{}9}}%
      \put(3399,0){\makebox(0,0){\strut{}10}}%
    }%
    \gplgaddtomacro\gplfronttext{%
      \csname LTb\endcsname
      \put(-1001,1809){\rotatebox{-270}{\makebox(0,0){\strut{}Average gap to company's solution}}}%
      \put(1699,-330){\makebox(0,0){\strut{}$K$}}%
      \csname LTb\endcsname
      \put(2544,2627){\makebox(0,0)[r]{\strut{}600s}}%
      \csname LTb\endcsname
      \put(2544,2407){\makebox(0,0)[r]{\strut{}1200s}}%
      \csname LTb\endcsname
      \put(2544,2187){\makebox(0,0)[r]{\strut{}1800s}}%
      \csname LTb\endcsname
      \put(2544,1967){\makebox(0,0)[r]{\strut{}2400s}}%
      \csname LTb\endcsname
      \put(2544,1747){\makebox(0,0)[r]{\strut{}3000s}}%
      \csname LTb\endcsname
      \put(2544,1527){\makebox(0,0)[r]{\strut{}3600s}}%
    }%
    \gplbacktext
    \put(0,0){\includegraphics{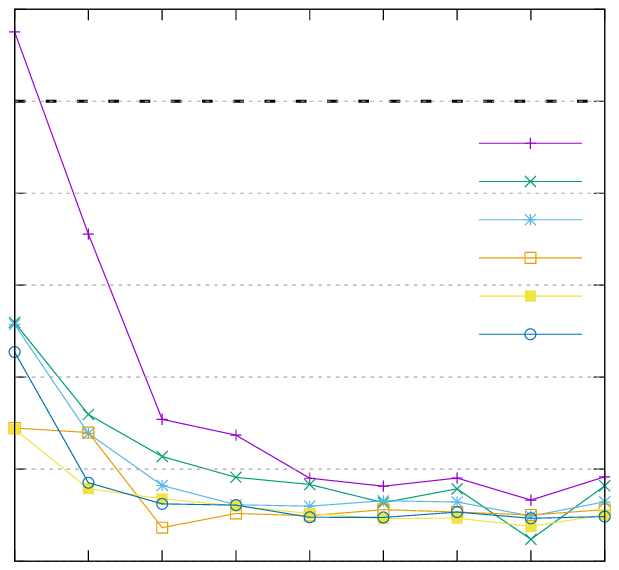}}%
    \gplfronttext
  \end{picture}%
\endgroup

%% file: abkm2020fig7b.tex
\begingroup
  \makeatletter
  \providecommand\color[2][]{%
    \GenericError{(gnuplot) \space\space\space\@spaces}{%
      Package color not loaded in conjunction with
      terminal option `colourtext'%
    }{See the gnuplot documentation for explanation.%
    }{Either use 'blacktext' in gnuplot or load the package
      color.sty in LaTeX.}%
    \renewcommand\color[2][]{}%
  }%
  \providecommand\includegraphics[2][]{%
    \GenericError{(gnuplot) \space\space\space\@spaces}{%
      Package graphicx or graphics not loaded%
    }{See the gnuplot documentation for explanation.%
    }{The gnuplot epslatex terminal needs graphicx.sty or graphics.sty.}%
    \renewcommand\includegraphics[2][]{}%
  }%
  \providecommand\rotatebox[2]{#2}%
  \@ifundefined{ifGPcolor}{%
    \newif\ifGPcolor
    \GPcolorfalse
  }{}%
  \@ifundefined{ifGPblacktext}{%
    \newif\ifGPblacktext
    \GPblacktexttrue
  }{}%
  \let\gplgaddtomacro\g@addto@macro
  \gdef\gplbacktext{}%
  \gdef\gplfronttext{}%
  \makeatother
  \ifGPblacktext
    \def\colorrgb#1{}%
    \def\colorgray#1{}%
  \else
    \ifGPcolor
      \def\colorrgb#1{\color[rgb]{#1}}%
      \def\colorgray#1{\color[gray]{#1}}%
      \expandafter\def\csname LTw\endcsname{\color{white}}%
      \expandafter\def\csname LTb\endcsname{\color{black}}%
      \expandafter\def\csname LTa\endcsname{\color{black}}%
      \expandafter\def\csname LT0\endcsname{\color[rgb]{1,0,0}}%
      \expandafter\def\csname LT1\endcsname{\color[rgb]{0,1,0}}%
      \expandafter\def\csname LT2\endcsname{\color[rgb]{0,0,1}}%
      \expandafter\def\csname LT3\endcsname{\color[rgb]{1,0,1}}%
      \expandafter\def\csname LT4\endcsname{\color[rgb]{0,1,1}}%
      \expandafter\def\csname LT5\endcsname{\color[rgb]{1,1,0}}%
      \expandafter\def\csname LT6\endcsname{\color[rgb]{0,0,0}}%
      \expandafter\def\csname LT7\endcsname{\color[rgb]{1,0.3,0}}%
      \expandafter\def\csname LT8\endcsname{\color[rgb]{0.5,0.5,0.5}}%
    \else
      \def\colorrgb#1{\color{black}}%
      \def\colorgray#1{\color[gray]{#1}}%
      \expandafter\def\csname LTw\endcsname{\color{white}}%
      \expandafter\def\csname LTb\endcsname{\color{black}}%
      \expandafter\def\csname LTa\endcsname{\color{black}}%
      \expandafter\def\csname LT0\endcsname{\color{black}}%
      \expandafter\def\csname LT1\endcsname{\color{black}}%
      \expandafter\def\csname LT2\endcsname{\color{black}}%
      \expandafter\def\csname LT3\endcsname{\color{black}}%
      \expandafter\def\csname LT4\endcsname{\color{black}}%
      \expandafter\def\csname LT5\endcsname{\color{black}}%
      \expandafter\def\csname LT6\endcsname{\color{black}}%
      \expandafter\def\csname LT7\endcsname{\color{black}}%
      \expandafter\def\csname LT8\endcsname{\color{black}}%
    \fi
  \fi
    \setlength{\unitlength}{0.0500bp}%
    \ifx\gptboxheight\undefined%
      \newlength{\gptboxheight}%
      \newlength{\gptboxwidth}%
      \newsavebox{\gptboxtext}%
    \fi%
    \setlength{\fboxrule}{0.5pt}%
    \setlength{\fboxsep}{1pt}%
\begin{picture}(3400.00,3400.00)%
    \gplgaddtomacro\gplbacktext{%
      \csname LTb\endcsname
      \put(-132,2922){\makebox(0,0)[r]{\strut{} }}%
      \put(0,0){\makebox(0,0){\strut{}2}}%
      \put(425,0){\makebox(0,0){\strut{}3}}%
      \put(850,0){\makebox(0,0){\strut{}4}}%
      \put(1275,0){\makebox(0,0){\strut{}5}}%
      \put(1700,0){\makebox(0,0){\strut{}6}}%
      \put(2124,0){\makebox(0,0){\strut{}7}}%
      \put(2549,0){\makebox(0,0){\strut{}8}}%
      \put(2974,0){\makebox(0,0){\strut{}9}}%
      \put(3399,0){\makebox(0,0){\strut{}10}}%
    }%
    \gplgaddtomacro\gplfronttext{%
      \csname LTb\endcsname
      \put(-473,1809){\rotatebox{-270}{\makebox(0,0){\strut{} }}}%
      \put(1699,-330){\makebox(0,0){\strut{}$K$}}%
      \csname LTb\endcsname
      \put(2544,2627){\makebox(0,0)[r]{\strut{}600s}}%
      \csname LTb\endcsname
      \put(2544,2407){\makebox(0,0)[r]{\strut{}1200s}}%
      \csname LTb\endcsname
      \put(2544,2187){\makebox(0,0)[r]{\strut{}1800s}}%
      \csname LTb\endcsname
      \put(2544,1967){\makebox(0,0)[r]{\strut{}2400s}}%
      \csname LTb\endcsname
      \put(2544,1747){\makebox(0,0)[r]{\strut{}3000s}}%
      \csname LTb\endcsname
      \put(2544,1527){\makebox(0,0)[r]{\strut{}3600s}}%
    }%
    \gplbacktext
    \put(0,0){\includegraphics{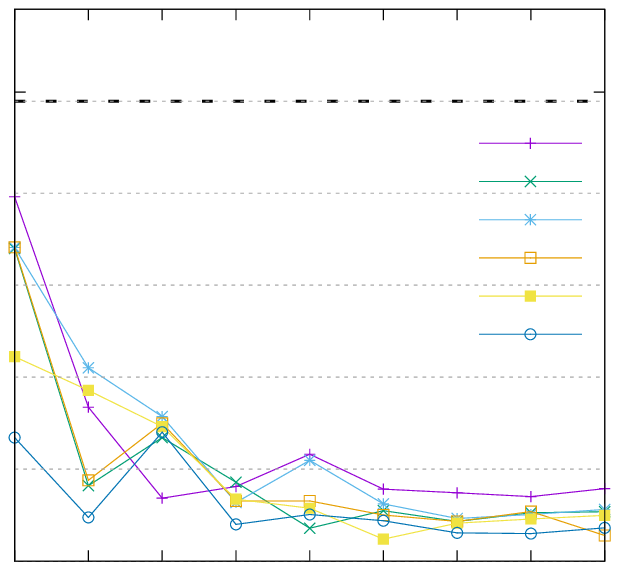}}%
    \gplfronttext
  \end{picture}%
\endgroup

%% file: abkm2020fig8.tex
\begingroup
  \makeatletter
  \providecommand\color[2][]{%
    \GenericError{(gnuplot) \space\space\space\@spaces}{%
      Package color not loaded in conjunction with
      terminal option `colourtext'%
    }{See the gnuplot documentation for explanation.%
    }{Either use 'blacktext' in gnuplot or load the package
      color.sty in LaTeX.}%
    \renewcommand\color[2][]{}%
  }%
  \providecommand\includegraphics[2][]{%
    \GenericError{(gnuplot) \space\space\space\@spaces}{%
      Package graphicx or graphics not loaded%
    }{See the gnuplot documentation for explanation.%
    }{The gnuplot epslatex terminal needs graphicx.sty or graphics.sty.}%
    \renewcommand\includegraphics[2][]{}%
  }%
  \providecommand\rotatebox[2]{#2}%
  \@ifundefined{ifGPcolor}{%
    \newif\ifGPcolor
    \GPcolorfalse
  }{}%
  \@ifundefined{ifGPblacktext}{%
    \newif\ifGPblacktext
    \GPblacktexttrue
  }{}%
  \let\gplgaddtomacro\g@addto@macro
  \gdef\gplbacktext{}%
  \gdef\gplfronttext{}%
  \makeatother
  \ifGPblacktext
    \def\colorrgb#1{}%
    \def\colorgray#1{}%
  \else
    \ifGPcolor
      \def\colorrgb#1{\color[rgb]{#1}}%
      \def\colorgray#1{\color[gray]{#1}}%
      \expandafter\def\csname LTw\endcsname{\color{white}}%
      \expandafter\def\csname LTb\endcsname{\color{black}}%
      \expandafter\def\csname LTa\endcsname{\color{black}}%
      \expandafter\def\csname LT0\endcsname{\color[rgb]{1,0,0}}%
      \expandafter\def\csname LT1\endcsname{\color[rgb]{0,1,0}}%
      \expandafter\def\csname LT2\endcsname{\color[rgb]{0,0,1}}%
      \expandafter\def\csname LT3\endcsname{\color[rgb]{1,0,1}}%
      \expandafter\def\csname LT4\endcsname{\color[rgb]{0,1,1}}%
      \expandafter\def\csname LT5\endcsname{\color[rgb]{1,1,0}}%
      \expandafter\def\csname LT6\endcsname{\color[rgb]{0,0,0}}%
      \expandafter\def\csname LT7\endcsname{\color[rgb]{1,0.3,0}}%
      \expandafter\def\csname LT8\endcsname{\color[rgb]{0.5,0.5,0.5}}%
    \else
      \def\colorrgb#1{\color{black}}%
      \def\colorgray#1{\color[gray]{#1}}%
      \expandafter\def\csname LTw\endcsname{\color{white}}%
      \expandafter\def\csname LTb\endcsname{\color{black}}%
      \expandafter\def\csname LTa\endcsname{\color{black}}%
      \expandafter\def\csname LT0\endcsname{\color{black}}%
      \expandafter\def\csname LT1\endcsname{\color{black}}%
      \expandafter\def\csname LT2\endcsname{\color{black}}%
      \expandafter\def\csname LT3\endcsname{\color{black}}%
      \expandafter\def\csname LT4\endcsname{\color{black}}%
      \expandafter\def\csname LT5\endcsname{\color{black}}%
      \expandafter\def\csname LT6\endcsname{\color{black}}%
      \expandafter\def\csname LT7\endcsname{\color{black}}%
      \expandafter\def\csname LT8\endcsname{\color{black}}%
    \fi
  \fi
    \setlength{\unitlength}{0.0500bp}%
    \ifx\gptboxheight\undefined%
      \newlength{\gptboxheight}%
      \newlength{\gptboxwidth}%
      \newsavebox{\gptboxtext}%
    \fi%
    \setlength{\fboxrule}{0.5pt}%
    \setlength{\fboxsep}{1pt}%
\begin{picture}(7200.00,7200.00)%
    \gplgaddtomacro\gplbacktext{%
      \csname LTb\endcsname
      \put(-132,220){\makebox(0,0)[r]{\strut{}$0\%$}}%
      \put(-132,854){\makebox(0,0)[r]{\strut{}$10\%$}}%
      \put(-132,1489){\makebox(0,0)[r]{\strut{}$20\%$}}%
      \put(-132,2123){\makebox(0,0)[r]{\strut{}$30\%$}}%
      \put(-132,2758){\makebox(0,0)[r]{\strut{}$40\%$}}%
      \put(-132,3392){\makebox(0,0)[r]{\strut{}$50\%$}}%
      \put(-132,4027){\makebox(0,0)[r]{\strut{}$60\%$}}%
      \put(-132,4661){\makebox(0,0)[r]{\strut{}$70\%$}}%
      \put(-132,5296){\makebox(0,0)[r]{\strut{}$80\%$}}%
      \put(-132,5930){\makebox(0,0)[r]{\strut{}$90\%$}}%
      \put(-132,6565){\makebox(0,0)[r]{\strut{}$100\%$}}%
      \put(-132,7199){\makebox(0,0)[r]{\strut{}$110\%$}}%
      \put(450,0){\makebox(0,0){\strut{}P1}}%
      \put(1350,0){\makebox(0,0){\strut{}P2}}%
      \put(2250,0){\makebox(0,0){\strut{}P3}}%
      \put(3150,0){\makebox(0,0){\strut{}P4}}%
      \put(4049,0){\makebox(0,0){\strut{}P5}}%
      \put(4949,0){\makebox(0,0){\strut{}P6}}%
      \put(5849,0){\makebox(0,0){\strut{}P7}}%
      \put(6749,0){\makebox(0,0){\strut{}P8}}%
      \put(3869,6438){\makebox(0,0)[l]{\strut{}$\scriptstyle  271.62\%$}}%
      \put(4634,6438){\makebox(0,0)[l]{\strut{}$\scriptstyle  118.81\%$}}%
      \put(5669,6438){\makebox(0,0)[l]{\strut{}$\scriptstyle  646.76\%$}}%
      \put(6209,6691){\makebox(0,0)[l]{\strut{}$\scriptstyle  646.76\%$}}%
      \put(6479,6438){\makebox(0,0)[l]{\strut{}$\scriptstyle  171.02\%$}}%
    }%
    \gplgaddtomacro\gplfronttext{%
      \csname LTb\endcsname
      \put(1215,6137){\makebox(0,0)[r]{\strut{}CPLEX}}%
      \csname LTb\endcsname
      \put(1215,5917){\makebox(0,0)[r]{\strut{}S1 (K = 8)}}%
      \csname LTb\endcsname
      \put(1215,5697){\makebox(0,0)[r]{\strut{}S9 (K = 8)}}%
      \csname LTb\endcsname
      \put(2250,6565){\makebox(0,0){\strut{}Company's solution}}%
    }%
    \gplbacktext
    \put(0,0){\includegraphics{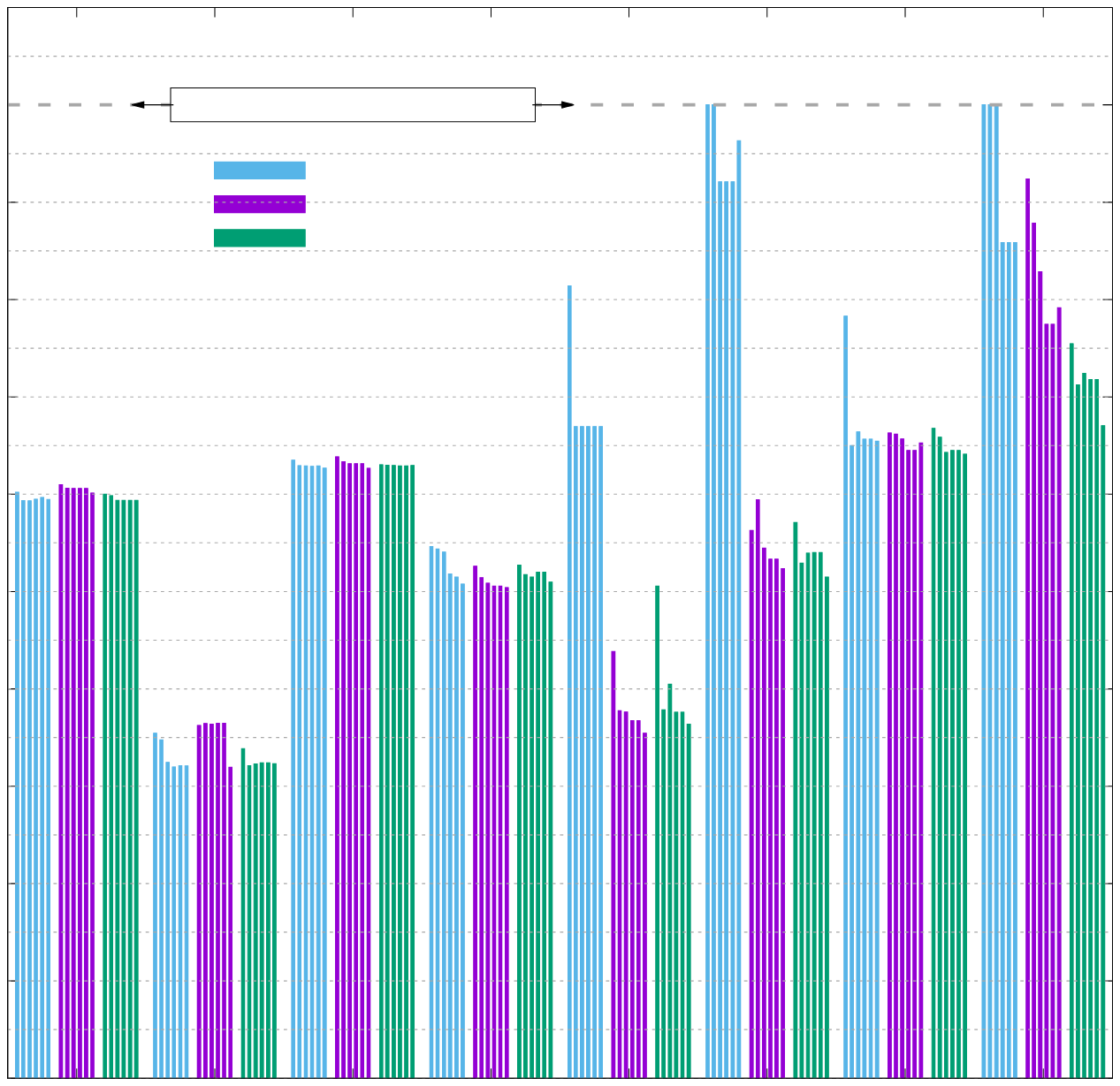}}%
    \gplfronttext
  \end{picture}%
\endgroup